\documentclass[12pt,reqno]{amsart}
\usepackage[margin=1.0in]{geometry} 
\usepackage{amscd, amssymb, amsmath,amsfonts}
\usepackage{mathrsfs}
\usepackage{graphicx}
\graphicspath{{./Figures/}}
\usepackage{subcaption,caption}
\usepackage[hidelinks]{hyperref}
\usepackage{thmtools,thm-restate}
\usepackage{tikz-cd}

\newtheorem{thm}{Theorem}[section]
\newtheorem{lem}[thm]{Lemma}

\theoremstyle{definition}
\newtheorem{rem}[thm]{Remark}
\newtheorem{defn}[thm]{Definition}
\newtheorem{ex}[thm]{Example}
\newtheorem{conv}[thm]{Convention}
\newcommand{\mvec}[1]{{\smash{\text{\boldmath{$#1$}}}}}
\renewcommand{\v}{\mvec{\lor}}
\newcommand{\n}{\mvec{\wedge}}
\newcommand{\tr}{\operatorname{tr}}
\newcommand{\qtr}{\operatorname{qtr}}
\newcommand{\End}{\operatorname{End}}
\newcommand{\floor}[1]{\left\lfloor #1 \right\rfloor}
\newcommand{\sltwo}{\mathfrak{sl}_2}
\newcommand{\C}{\mathbb{C}}
\newcommand{\K}{\mathbb{K}}

\newcommand{\Z}{\mathbb{Z}}
\newcommand{\R}{\mathcal{R}}
\newcommand{\bracks}[1]{\left\langle #1\right\rangle}
\DeclareMathOperator{\Hom}{Hom}
\title[Seifert-Torres Type Formulas from Quantum $\mathfrak{sl}_2$]
{Seifert-Torres Type Formulas for the Alexander Polynomial from Quantum $\mathfrak{sl}_2$}
\author{Matthew Harper}
\address{\parbox{0.9\linewidth}{Department of Mathematics, University of California, Riverside, CA 92521, USA \hfill \\
Department of Mathematics, The Ohio State University, 
		Columbus, OH 43210, USA}}
\email{mharper@ucr.edu}
\begin{document}
	\maketitle\thispagestyle{empty}
	\begin{abstract} \noindent
		We develop a diagrammatic calculus for representations of unrolled quantum $\mathfrak{sl}_2$ at a fourth root of unity. This allows us to prove Seifert-Torres type formulas for certain splice links using quantum algebraic methods, rather than topological methods. Other applications of this diagrammatic calculus given here are a skein relation for $n$-cabled double crossings and a simple proof that the quantum invariant associated with these representations determines the multivariable Alexander polynomial.
	\end{abstract}
\section{Introduction}
The Alexander polynomial is an invariant of knots first constructed by J. W. Alexander \cite{Alexander}. It was refined by Fox to an invariant of links, which we call the multivariable Alexander polynomial $\Delta$ \cite{Fox}. For an $n$-component link, $\Delta$ takes values in $\mathbb{Z}[t_1^\pm,\dots, t_n^\pm]$ and is defined up to a unit. The Conway Potential Function (CPF) $\nabla$ is a well-defined standard representative of the multivariable Alexander polynomial as a rational function in the following sense \cite{Conway}. For every oriented link $L$ with colored (ordered) components,
\begin{align}
\nabla(\mathcal{L})(t_1,\dots, t_n)\dot{=}
\begin{cases}
	\dfrac{\Delta(\mathcal{L})(t_1^2)}{t_1-t_1^{-1}}& n=1\\
	\Delta(\mathcal{L})(t_1^2,\dots, t_n^2) & n>1,
\end{cases}
\end{align}
with a representative $\Delta(\mathcal{L})$ chosen such that $\Delta(\mathcal{L})(t_1^{-2},\dots, t_n^{-2})=\pm\Delta(\mathcal{L})(t_1^2,\dots, t_n^2)$. Our use of $\dot{=}$ indicates that the sign of the representative must also be chosen correctly in order to produce the CPF, this choice is specific to the link. In this paper, we define the CPF as the unique function on links that satisfies the Jiang relations given here in Figure \ref{fig:rels} \cite{Jiang}.\\

Jun Murakami constructs $\nabla$ as a quantum invariant in \cite{MurakamiStateModel} as a Turaev-type state model \cite{Turaev88}. In \cite{Murakami}, he shows that the $R$-matrix used in his state model coincides with a (parameterized) representation of the $R$-matrix of (unrolled) restricted quantum $\sltwo$ at a primitive fourth root of unity.\\

We develop a diagrammatic calculus for these quantum group representations, which we then apply to prove several results for the CPF. 

\begin{restatable*}{thm}{sat}\label{thm:sat}
	Let $\mathcal{T}_1,\dots, \mathcal{T}_n$ be colored string links where $\mathcal{T}_i$ has $d_i$ components and the color of its $j$-th component is $c_{ij}\in\{t_{i,1},\dots, t_{i,m_i}\}$. Let $\mathcal{L}$ be a framed $n$-component link with zero linking matrix and $\mathcal{L}'$ the $(\mathcal{L},\mathcal{T}_1,\dots, \mathcal{T}_n)$-satellite link. Set $\tau_i=c_{i, 1} \cdots c_{i,d_i}$.  Then
	\begin{align*}
		\nabla(\mathcal{L}')(t_{1,1},\dots, t_{n,m_n})=\nabla(\mathcal{L})(\tau_1,\dots, \tau_n) \prod_{i=1}^n\left(\nabla(\widehat{\mathcal{T}_i})(t_{i, 1}, \dots, t_{i,m_i})(\tau_i-\tau_i^{-1})\right).
	\end{align*}
\end{restatable*} 

The link $\mathcal{L}'$ in the above theorem is a particular example of a splice link \cite[Chapter 1]{EN}. Let $\mathcal{S}_i$ be the link $\widehat{\mathcal{T}_i}\cup S^1$, where $S^1$ is an unknot homologous to the meridian of the torus naturally containing $\widehat{\mathcal{T}}_i$. Then $\mathcal{L}'$ is obtained by splicing each component of $\mathcal{L}$ with the unknot in the corresponding ${\mathcal{S}_i}$.\\

 For our purposes, it is easier to consider this from a diagrammatic perspective. We take a framed string link $\mathcal{L}_0$ with $n$ components obtained from a partial closure of a braid representative of $\mathcal{L}$. On each of the  components of $\mathcal{L}_0$, we take its $d_i$-parallel cabling with respect to its framing and denote the resulting tangle $\mathcal{L}_0^{\parallel \overline{d}}$. We obtain the $(\mathcal{L},\mathcal{T}_1,\dots, \mathcal{T}_n)$-satellite link by taking the closure of the composite, $\mathcal{L}_0^{\parallel \overline{d}}$ stacked on $\mathcal{T}_1,\dots, \mathcal{T}_n$. The coloring on the tangles induces a coloring on the resulting link.\\

Theorem \ref{thm:sat} is an extension of the Seifert-Torres formula \cite{Seifert,Torres56} and is a corollary of Cimasoni's formula for a splice link  \cite[Theorem 3.1]{splice}. Although this theorem is not new, we use techniques from quantum algebra rather than topology in our proofs.\\

A key lemma in the theorem's proof is a naturality result between the action of a string link with zero linking matrix on a representation with product colors and its cabling acting on a tensor product. This is most easily presented as a commutative diagram for a framed long knot diagram (1-tangle) ${\mathcal{K}_0}$. The diagram commutes if and only if ${\mathcal{K}_0}$ has zero writhe. 

\begin{figure}[h!]
	\centering
	\begin{tikzcd}
		V( t_1\cdots t_n) \arrow[r, "\bracks{\mathcal{K}_0}"] \arrow[d,,swap, hook]
		& 	V(t_1\cdots t_n) \arrow[d,hook ] \\
		\bigotimes_{i=1}^{n} V(t_i) \arrow[r, "\bracks{\mathcal{K}_0^{\parallel n}}"]
		& \bigotimes_{i=1}^{n} V(t_i)
	\end{tikzcd}
\end{figure}
For links with arbitrary linking matrix, we prove the following theorem. 
\begin{restatable*}{thm}{link}\label{thm:link}
	Let $\mathcal{L}$ be a framed $n$-component link with linking matrix $(l_{ij})$. We replace each of the $i$ components of $\mathcal{L}$ with a $d_i$-fold parallel-cable along its framing and denote the resulting link by $\mathcal{L}^{\parallel \overline{d}}$. Let $t_{j,k}$ be the color of the $k$-th cable of the $j$-th component. Set $\tau_j=t_{j,1}\cdots t_{j,d_i}$. Then 
	\begin{align*}
		\nabla(\mathcal{L}^{\parallel \overline{d}})(t_{1,1},\dots, t_{n,d_n})=&\nabla(\mathcal{L})(\tau_1,\dots,\tau_n)\prod_{i=1}^n \left(\prod_{j=1}^n\tau_j^{l_{ij}}-\prod_{j=1}^n\tau_j^{-l_{ij}}\right)^{d_i-1}.
	\end{align*}
\end{restatable*}
Theorem \ref{thm:link} is a special case of the formula for $\nabla_{L'\setminus K_n}$ in \cite[Corollary 3.4]{splice}, as $\mathcal{L}^{\parallel \overline{d}}$ is obtained from splicing $\mathcal{L}$ with links consisting of $d_i$ copies of $(1,l_{ii})$-torus knots union an unknot. The splicing operation assumes the components involved have standard longitudes and meridians. Therefore, taking a $(1,q)$-cabling, as in \cite{splice}, is equivalent to taking the framed cabling along a $q$-framed component, as in the context of the above theorem.\\

    We also apply the $\sltwo$ diagrammatic calculus to give a much simpler and more concise proof of Murakami's theorem \cite{MurakamiStateModel} by verifying a set of axioms for the CPF discovered by Jiang \cite{Jiang}, built on earlier work by Hartley \cite{Hartley}. These axioms are given in Figure \ref{fig:rels}. 
    Jiang's approach to characterizing the CPF is entirely topological, and avoids the representation theoretic tools employed in Murakami's work.
      
   \begin{thm}[\cite{MurakamiStateModel}]\label{thm}
   	The link invariant determined by the representations $V(t)$ of $\overline{U}_\zeta^H(\sltwo)$ is the Conway Potential Function.
   \end{thm}

We also show that Jiang relation (\ref{eq:II}) generalizes to cabled strings. As in Figure \ref{fig:Rpar}, let $R_{t_1,\dots, t_n}^{\parallel n}\in\End\left(\bigotimes_{i=1}^n V(t_i)\otimes \bigotimes_{i=1}^n V(t_i)\right)$ denote the linear map associated to a crossing between two sets of $n$-cabled strands.
\begin{restatable*}{thm}{cableskein}
	The following skein relation holds for cabled strands:
	\begin{equation*}
		\left(R_{t_1,\dots, t_n}^{\parallel n}\right)^2+\left(R_{t_1,\dots, t_n}^{\parallel n}\right)^{-2}=
		\left(t_1^2\cdots t_n^2+t_1^{-2}\cdots t_n^{-2}\right)
		id_{\bigotimes V(t_i)\otimes \bigotimes V(t_i)}.	\end{equation*}
\end{restatable*}

\subsection{Related Work}
Akutsu, Deguchi, and Viro study diagrammatic calculi for quantum $\mathfrak{sl}_2$ similar to the one in this paper, but with different conventions \cite{DeguchiAkutsu,Viro}. A minor difference is their use of the additive weight convention in representations, rather than multiplicative, so computations are not written in terms of the explicit polynomial variables. \\

A more significant difference in our work is in regards to forks/trivalent vertices. Although the calculus described here defines an invariant of labeled trivalent graphs, our focus is on the link invariant. Viro's sink (source) vertex is what we call a positive (negative) fork. Forks considered in this paper are normalized to satisfy orthonormality conditions as in the paper by Deguchi and Akutsu, unlike Viro's convention.  As pointed out in \cite[Section 3]{Viro}, a disadvantage of this method is that at each vertex one of the adjacent edges is distinguished. However, since we work with string links and braids, which have a natural upward orientation, it is more convenient for us to use the orthonormalized version. We have also worked out associator moves for forks and their interaction with normalized $R$-matrices. The author is not aware of these computations in the present literature.\\

We also note that Viro uses Turaev's axioms \cite{Turaev86} to show that the $\overline{U}_\zeta^H(\mathfrak{sl}_2)$ invariant is the CPF. Here we use Jiang's skein relations \cite{Jiang} for the CPF which are diagrammatic and are more easily verified with the calculus we develop.

\subsection{Future Work} A technical barrier to extending our results to the full Seifert-Torres formula, is the indecomposability of $V(t)\otimes V(t)^*$, which we avoid by considering string links and assuming non-degeneracy between all pairs of parameters. These representations are not compatible with the forks used in the current diagrammatic calculus. If they can be incorporated,
we expect our proof of {Theorem~\ref{thm:sat}} can be applied to more general
splice links and that we can extend
Rolfsen's formula for the Alexander polynomial of Whitehead doubles
\cite{Rolfsen}. 

\subsection{Acknowledgments}
I thank Sergei Chmutov for telling me about Jiang's work on the multivariable Alexander polynomial. I am also grateful to Thomas Kerler for helpful comments. I also thank an anonymous referee for pointing out \cite{splice}. I thank the NSF for partial support through the grant \textbf{NSF-RTG \#DMS-1547357}.

\section{Quantum Algebra Background}
\begin{figure}[h!]
	\centering
	\includegraphics[scale=1]{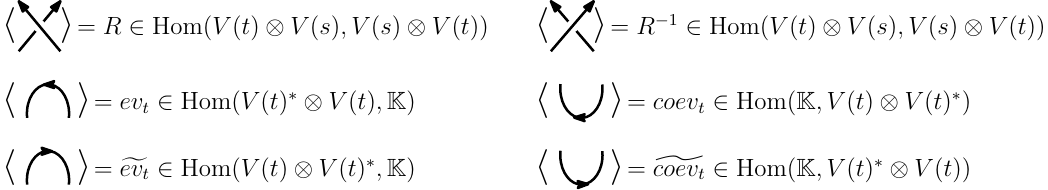}
	\caption{Linear maps assigned to crossings, cups, and caps in a diagram.}\label{fig:RT}
\end{figure}

In this section we recall key facts about quantum invariants of links from Hopf algebras, recall the unrolled restricted quantum group $\overline{U}_\zeta^H(\sltwo)$ and its representations, and establish conventions for the paper. Given a representation of a ribbon Hopf algebra, such as  $\overline{U}_\zeta^H(\sltwo)$, we may assign intertwiner maps to elementary tangles via the Reshetikhin-Turaev functor, see Figure \ref{fig:RT}. Most of the content in this section can be found in \cite[Chapter 4]{Ohtsuki} and the reader is referred there for additional details. Many of our conventions follow Ohtsuki's; however, our conventions differ in that an upward strand is associated to $id_V$; not $id_{V^*}$, and our ground ring for the quantum group is $\Z(\zeta)$; not $\Z(\zeta^{1/2})$, where $\zeta$ is a primitive fourth root of unity.\\

A Hopf algebra over a field $\K$ is an algebra $A$ equipped with homomorphisms $\Delta:A\to A\otimes A $ and $\epsilon :A\to \mathbb{K}
$, and an anti-homomorphism $S:A\to A
$ such that 
\begin{align}
	(id\otimes \Delta)\circ\Delta&=\Delta^2=(\Delta\otimes id)\circ\Delta,\\
	(id\otimes \epsilon)\circ\Delta&=id=(\epsilon\otimes id)\circ\Delta,
	\\
	m\circ(id\otimes S)\circ \Delta&=i \circ \epsilon=m\circ(S\otimes id)\circ \Delta,
\end{align}
where $m:A\otimes A\to A$ is multiplication on $A$ and $i:\K\to A$ is the unit map. The maps $\Delta,~\epsilon,$ and $S$ are called the coproduct, counit, and antipode, respectively. The counit defines the trivial representation of $A$. The coproduct and antipode define the action of $A$ on tensor product and dual representations, receptively.\\
 
\begin{conv}
For the remainder of the paper $\Delta$ will denote the coproduct on a Hopf algebra, not the Alexander polynomial. Moreover, we use a calligraphic font to identify topological objects to provide additional context for these two maps.
\end{conv}

We will focus on the Hopf algebra $\overline{U}_\zeta^H(\sltwo)$ the unrolled restricted quantum group at a primitive fourth root of unity $\zeta$. We will assume $\zeta=e^{i\pi/2}$, but its conjugate may be used just as well.

\begin{defn}
	Let $\overline{U}_\zeta^H(\sltwo)$ be the $\mathbb{Z}(\zeta)$-algebra generated by $E,F,H,K,K^{-1}$ with relations
	\begin{align}
		KK^{-1}=K^{-1}K=1, && KE=-EK, && KF=-FK, && [E,F]=\frac{K-K^{-1}}{\zeta-\zeta^{-1}},\\
		HK=KH, && [H,E]=2E,&&[H,F]=-2F,&& E^2=F^2=0.		
	\end{align}
\end{defn}
The Hopf algebra structure on $\overline{U}_\zeta^H(\sltwo)$ is given by
\begin{align}
	\Delta(E)&=E\otimes K+1\otimes E, & S(E)&=-EK^{-1},&\epsilon(E)&=0,\label{eq:HopfE}\\
	\Delta(F)&=F\otimes 1+ K^{-1}\otimes F, & S(F)&=-KF,&\epsilon(F)&=0,\label{eq:HopfF}\\
	\Delta(H)&=H\otimes 1+1\otimes H, & S(H)&=-H,&\epsilon(H)&=0,\label{eq:HopfH}\\
	\Delta(K)&=K\otimes K,&S(K)&=K^{-1},&\epsilon(K)&=1\label{eq:HopfK}.
\end{align}

\begin{defn}
	A Hopf algebra $A$ is said to be quasi-triangular if there exists an invertible element $\R= \sum_i \alpha_i\otimes\beta_i\in A\otimes A$ called the universal $R$-matrix such that
	\begin{align}
		&\Delta^{\text{op}}(a)=\R \Delta(a) \R^{-1}~ \text{ for any }a\in A,
		\\ &(\Delta \otimes id )(\R)=\R_{13}\R_{23}, \label{eq:R1323}
		\\ &(id\otimes \Delta)(\R)=\R_{13}\R_{12}. \label{eq:R1312}
	\end{align}
	By $\R_{ij}$ we mean  $\R_{12}=\R\otimes 1$, $\R_{23}=1\otimes \R$, and $\R_{13}=\sum_i \alpha_i\otimes1\otimes\beta_i$. We write $\Delta^{\text{op}}$ to mean $\tau\circ \Delta$, where $\tau$ is the tensor swap $\tau(x\otimes y)=y\otimes x$. 
\end{defn}

As a consequence of (\ref{eq:R1323}) and (\ref{eq:R1312}), universal $R$-matrices satisfy the Yang-Baxter equation
\begin{align}
	\R_{12}\R_{13}\R_{23}=\R_{23}\R_{13}\R_{12}.
\end{align}
This relation connects Hopf algebras and braid groups, as it implements the third Reidemeister move.\\

For ease of exposition, we postpone the discussion of the universal $R$-matrix of $\overline{U}_\zeta^H(\sltwo)$, which has a simpler description on representations.

\begin{defn}
	Given a quasi-triangular Hopf algebra, the Drinfeld element, defined in	terms of the $R$-matrix, is given by $u=\sum_i S(\beta_i)\alpha_i$.
\end{defn}

\begin{defn}
	A quasi-triangular Hopf algebra $A$ is said to be \textit{ribbon} if it contains a central element $v\in A$ such that
	\begin{align}
		v^2=S(u)u, &&\Delta(v)=(v\otimes v)\cdot (\R_{21}\R)^{-1}, &&S(v)=v, && \epsilon(v)=1.
	\end{align}
\end{defn}
Such a $v$ is called a \textit{ribbon element} and we call $uv^{-1}$ the \textit{pivotal element}. For $\overline{U}_\zeta^H(\sltwo)$, the pivotal element is $K^{-1}$.\newline

Given a finite-dimensional representation $(V,\rho_V)$ of a ribbon Hopf algebra $A$ and an oriented link diagram $\mathcal{L}$, the Reshetikhin-Turaev functor assigns a morphism in $\operatorname{Rep}A$ to $\mathcal{L}$ which is determined by composing the morphisms assigned to elementary tangles which comprise $\mathcal{L}$. This assignment is given in Figure \ref{fig:RT} and we use the convention that an upward pointing strand is the identity on $V$ and a downward pointing strand is the identity on $V^*$. If $\mathcal{L}$ has multiple components, each component of $\mathcal{L}$ may be ``colored'' by a different representation. Suppose that $(W,\rho_W)$ is also a representation of $A$, then $R$ denotes the intertwiner $\tau\circ(\rho_V\otimes \rho_W)(\R)$. The image of $\mathcal{L}$ under the Reshetikhin-Turaev functor is an invariant of the underlying link \cite{RT90,RT}, see also \cite[Theorem 4.7]{Ohtsuki}. This construction also produces an isotopy invariant of tangles.\\

Since the image of the functor is an isotopy invariant, evaluation and coevaluation satisfy the duality (zig-zag) relations
\begin{align}\label{eq:zigzag1}
	(id_{V}\otimes ev_V)(coev_V\otimes id_{V})&=id_{V}=(\widetilde{ev_V}\otimes id_{V} )(id_{V}\otimes \widetilde{coev_V}),\\
	\label{eq:zigzag2}
	(ev_V\otimes id_{V^*} )(id_{V^*}\otimes coev_V)&=id_{V^*}=	(id_{V^*}\otimes \widetilde{ev_V})(\widetilde{coev_V}\otimes id_{V^*}).
\end{align}
Given any basis $(e_i)$ of $V$ and corresponding dual basis $(e_i^*)$, the above maps are defined as 
\begin{align}
	ev_V(e_i^*\otimes e_j)=e_i^*(e_j),&& \widetilde{ev_V}(e_i\otimes e_j^*)= e_j^*\left(\rho(uv^{-1}) e_i\right),\\
	coev_V(1)=\sum_i e_i\otimes e_i^*, && \widetilde{coev_V}(1)= \sum_i e_i^*\otimes \rho(uv^{-1})^{-1} e_i,
\end{align}
and do not depend on the choice of basis.
\\

We consider invariants from $\overline{U}_\zeta^H(\sltwo)$ obtained from the following representations. 
\begin{defn}
	Fix $t\in\C^\times$. Let $(\rho_t,V(t))$ be a two-dimensional representation of $\overline{U}_\zeta^H(\sltwo)$, expressed in the standard basis $\langle v_0^t,Fv_0^t\rangle$ as
	\begin{align}
		\rho_t(F)=\begin{bmatrix}
			0&0\\1&0
		\end{bmatrix}, &&
		\rho_t(E)=\frac{t-t^{-1}}{\zeta-\zeta^{-1}}\begin{bmatrix}
			0&1\\0&0
		\end{bmatrix}, &&
		\rho_t(K)=\begin{bmatrix}
			t&0\\0&-t
		\end{bmatrix},
		&&
		\rho_t(H)=\begin{bmatrix}
			\alpha&0\\0&\alpha-2
		\end{bmatrix}
	\end{align}
	for some $\alpha$ such that $t=\zeta^{\alpha}=e^{i\pi\alpha/2}$. 
\end{defn}

We write $\floor{t}$ to mean $\dfrac{t-t^{-1}}{\zeta-\zeta^{-1}}$. With this notation, $\rho_t(E)=\begin{bmatrix}
	0&\floor{t}\\0&0
\end{bmatrix}$.

\begin{lem}
	 The representation $V(t)$ is irreducible if and only if $\floor{t}\neq 0.$
\end{lem}

\begin{lem}\label{lem:dual}
	The dual representation $V(t)^*$ is isomorphic to $V(-t^{-1})$ for all $t\in\C^\times$. 
\end{lem}

 With $v_{1}^{t}$ denoting $Fv_{0}^{t}$, the standard basis of $V(t)$ extends to a basis of $V(t)\otimes V(s)$ as \begin{align}\label{eq:basis}
	\langle v_0^t\otimes v_0^s,v_0^t\otimes v_1^s,v_1^t\otimes v_0^s,v_1^t\otimes v_1^s\rangle.
\end{align}

There is an action of the braid group on the category tensor-generated by the representations $\{V(t):t\in\C^\times \}$, making it into a braid groupoid. The braiding, which is defined in terms of $H$, can be normalized to only depend on $t=\zeta^\alpha$ and $s=\zeta^\beta$. Therefore, the precise values of $\alpha$ and $\beta$ will not be important for our discussion, as shown below.\\

The braid generators act via the $R$-matrix, $\label{eq:R}
R=\tau\circ\rho_t\otimes\rho_s
\left(\zeta^{H\otimes H/2}(1\otimes 1+(\zeta-\zeta^{-1}) E\otimes F)\right).
$
Here, $\tau$ is the tensor swap and $\zeta^{H\otimes H/2}(w_1\otimes w_2)=\zeta^{\lambda_1\lambda_2/2}(w_1\otimes w_2)$
for  $H$-weight vectors $w_1$ and $w_2$ of weight $\lambda_1$ and $\lambda_2$, respectively. In the standard tensor product basis, 

\allowdisplaybreaks
\begin{align}
	R&=\begin{bmatrix}
		1 & 0&0&0\\
		0& 0&1&0\\
		0&1&0&0\\
		0&0&0&1
	\end{bmatrix}\begin{bmatrix}
		\zeta^{\alpha\beta/2} & 0&0&0\\
		0&\zeta^{\alpha(\beta-2)/2}&0&0\\
		0& 0&\zeta^{(\alpha-2)\beta/2}&0\\
		0&0&0&\zeta^{(\alpha-2)(\beta-2)/2}
	\end{bmatrix}\begin{bmatrix}
		1 & 0&0&0\\
		0& 1&(t-t^{-1})&0\\
		0&0&1&0\\
		0&0&0&1
	\end{bmatrix}\\
	&=\zeta^{\alpha(\beta-2)/2}\begin{bmatrix}
		\zeta^{\alpha} & 0&0&0\\
		0& 0&\zeta^{\alpha-\beta}&0\\
		0&1&(t-t^{-1})&0\\
		0&0&0&-\zeta^{-\beta}
	\end{bmatrix}=\zeta^{\alpha(\beta-2)/2}\begin{bmatrix}
		t & 0&0&0\\
		0& 0&ts^{-1}&0\\
		0&1&t-t^{-1}&0\\
		0&0&0&-s^{-1}
	\end{bmatrix}. \label{eqn:Rts}
\end{align} 

\begin{conv}\label{conv:R}
	We define $R_{ts}:V(t)\otimes V(s)\to V(s)\otimes V(t)$ to be $\zeta^{-\alpha(\beta-2)/2}R$, which depends only on $t$ and $s$. We assign $R_{ts}$ to be the image of a positive crossing under the Reshetikhin-Turaev functor rather than $R$, see Figure \ref{fig:RT}.
\end{conv}
 Note $R_{ts}$ is an intertwiner and this rescaling does not change qualitative properties of the invariant. However, some naturality is lost as noted in Remark \ref{rem:notnatural}.\\

The following is a mild extension of a theorem of \cite[Lemma A.16]{Ohtsuki}.

\begin{thm}\label{thm:decomp}
	If $(ts)^2\neq 1$ or $t^2=s^2=1$, then $V(t)\otimes V(s)\cong V(ts)\oplus V(-ts)$. Otherwise, $V(t)\otimes V(s)$ is indecomposable and contains a 1-dimensional subrepresentation.
\end{thm}

\begin{rem}
	There are two isomorphism classes of such indecomposable tensor products. We denote them by $P(1)$ and $P(-1)$, as they are the projective covers of $V(1)$ and $V(-1)$, respectively. Moreover, $P(1)\cong{V(t)\otimes V(t^{-1})}$ and $P(-1)\cong V(t)\otimes V(-t^{-1})$ for all $t^2\neq 1 $.
\end{rem}

We say that the pair of parameters $(t,s)$ is \textit{degenerate} if the tensor product $V(t)\otimes V(s)$ is indecomposable. 

\begin{conv}\label{conv:nondeg}
	For the remainder of this paper, we assume all representations are irreducible and all pairs of parameters are non-degenerate unless noted otherwise. 
\end{conv}

\section{Tangles, Colorings, and Invariants from $\overline{U}_\zeta^H(\sltwo)$}\label{sec:color}
In this section we establish our convention for coloring tangles and describe how to compute a nontrivial invariant from the above described representations of $\overline{U}_\zeta^H(\sltwo)$. We prove in Section \ref{sec:CPFproof} that this invariant is the CPF.\newline
 
 A tangle $\mathcal{T}$ with $k$ lower boundary points and $l$ upper boundary points is called a $(k,l)$-tangle. If $k=l$, then we call $\mathcal{T}$ a $k$-tangle. Most tangles discussed in this paper are either links or (not necessarily pure) string links. Recall that a \textit{string link} is any embedding of $n$ oriented intervals into the cylinder $D^2\times [0,1]$ such that the initial and terminal points are mapped to $ D^2\times \{0\}$ and $ D^2\times \{1\}$, respectively. A string link with $n$ components is an example of an $n$-tangle.\\
 
 A \textit{coloring} of a tangle $\mathcal{T}$ is a surjection $c$ from its $n$ components to $\{t_1,\dots, t_m\}$ with the color of the $i$-th component of $\mathcal{T}$ is denoted $c_i$. The set $\{t_1,\dots, t_m\}$ is called the \textit{pallete} for the coloring. Let $\overline{c}=(c_1,\dots, c_k)$ and $\overline{c}'=(c'_1,\dots c'_l)$ be sequences of colors determined from the lower and upper boundaries of $\mathcal{T}$. If $\overline{c}=\overline{c}'$, then the components of $\widehat{\mathcal{T}}$ have a well-defined coloring. For string links, if $n$ is greater than the number of components of $\widehat{\mathcal{T}}$, $m<n$ and $c_i=c_j$ for some $i\neq j$. We usually choose a coloring so that all components of either $\mathcal{T}$ or $\widehat{\mathcal{T}}$ have a unique color. \newline 
 
 Composites of colored tangles $\mathcal{T}_1$ and $\mathcal{T}_2$ are also well-defined, provided the colors $\overline{c}'_1$ along the upper boundary of $\mathcal{T}_1$ agree with the colors $\overline{c}_2$ along the lower boundary of $\mathcal{T}_2$. Therefore, colored tangles can be written as the composite of colored elementary tangles.\\ 
 
 Each colored oriented $(k,l)$-tangle $\mathcal{T}$ determines a morphism in $\operatorname{Rep}\overline{U}_\zeta^H(\sltwo)$. For such a tangle there are induced orientations on its upper and lower boundaries, which we denote by $\epsilon\in \{\pm 1\}^k$ and $\epsilon'\in \{\pm 1\}^l$. Boundary points with an upward pointing orientation are assigned ``$+1$'' and are otherwise assigned ``$-1$.'' For $\overline{c}$ and $\overline{c}'$ as above, $\mathcal{T}$ determines a morphism $\bracks{\mathcal{T}}\in\Hom \left(\bigotimes_{i=1}^k  V(\epsilon_i {c_i}^{\epsilon_i}), \bigotimes_{i=1}^l V(\epsilon'_i{c'_i}^{\epsilon'_i})  \right)$ via the Reshetikhin-Turaev functor. Here, we have made the identification $V(t)^*\cong V(-t^{-1})$ from Lemma \ref{lem:dual}. The map $\bracks{\mathcal{T}}$ is defined over $\mathbb{K}:=\mathbb{Z}[\zeta, t_{1}^\pm,\dots, t_m^{\pm}]$ and is given by the composition of maps assigned to elementary oriented colored tangles comprising $\mathcal{T}$. \newline 
 
  If $\bracks{\mathcal{T}}\in\End\left(\bigotimes_{i=1}^k  V(\epsilon_i {c_i}^{\epsilon_i})\right)$, then its canonical trace $\tr:\End\left(\bigotimes_{i=1}^k  V(\epsilon_i {c_i}^{\epsilon_i})\right)\to\K$ and its right partial canonical trace $\tr_R:\End\left(\bigotimes_{i=1}^k  V(\epsilon_i {c_i}^{\epsilon_i})\right)\to\End\left(\bigotimes_{i=1}^{k-1}  V(\epsilon_i {c_i}^{\epsilon_i})\right)$ are well-defined. For simplicity, assume that $\mathcal{T}$ is a string link so that $\epsilon_i=1$ for all $i$. Also set $V'=\bigotimes_{i=1}^{k-1}  V(c_i)$.
	\begin{defn}
		 The \textit{right partial quantum trace} of $\bracks{\mathcal{T}}$, written $\qtr_R(\bracks{\mathcal{T}})$, is an endomorphism of $V'$ given by 
		\begin{align}
		\tr_R\left((id_{V'}\otimes \rho_{c_k}(K^{-1})) \bracks{\mathcal{T}}\right)=(id_{V'}\otimes \widetilde{ev_{V(c_k)}})\circ (\bracks{\mathcal{T}}\otimes id_{V(c_k)^*})\circ (id_{V'}\otimes coev_{V(c_k)}).
		\end{align}
		It is the image of the right strand diagrammatic closure of $\mathcal{T}$ under the Reshetikhin-Turaev functor. 
		The \textit{left partial quantum trace} $\qtr_L(\bracks{\mathcal{T}})$ is defined similarly.
	\end{defn}
	The operator obtained by taking $k-1$ successive right quantum partial traces $(\qtr_R)^{k-1}(\bracks{\mathcal{T}})$ acts by a scalar $\overline{\bracks{\mathcal{T}}}$ on $V(c_1)$ by the irreducibility assumption in Convention \ref{conv:nondeg}. Since $\tr(\rho_{c_1}(K^{-1}))=0$,
	\begin{align}
		\qtr_R^{k}(\bracks{\mathcal{T}})=\qtr_R((\qtr_R)^{k-1}(\bracks{\mathcal{T}}))=\tr(\rho_{c_1}(K^{-1})\overline{\bracks{\mathcal{T}}})=\overline{\bracks{\mathcal{T}}}\tr(\rho_{c_1}(K^{-1}))=0
	\end{align}
	and
	\begin{align}\label{eq:trace}
		\tr(\qtr_R^{k-1}(\bracks{\mathcal{T}}))=\tr\big(\overline{\bracks{\mathcal{T}}}\cdot id_{V(c_1)}\big)=\overline{\bracks{\mathcal{T}}}\dim\left(V(c_1)\right).
	\end{align}
	Thus, any closed tangle is assigned the zero morphism when colored by representations $V(t)$. However, if all but one strands are closed, then we obtain something nontrivial.\newline
	
	If a link $\mathcal{L}$ is given by the closure of some $k$-tangle $\mathcal{T}$, we would like $\tr_R^{k-1}(\bracks{T})$ to be an invariant of $\mathcal{L}$. Moreover, if $\mathcal{L}$ has multiple components, it should not depend on which component of $\mathcal{T}$ is left open. The additional steps needed to obtain an invariant of $\mathcal{L}$ from a tangle representative are given in \cite{GPT}, specifically in Section 6.3 where we take $N=2$. The following theorem summarizes their results for our purposes.
	
	\begin{thm}[{\cite{GPT}}]
		Let $\mathcal{L}$ be a link given by the closure of a $k$-component string link $\mathcal{T}$. Then \begin{align}
		\dfrac{\tr(\qtr_R^{k-1}(\bracks{\mathcal{T}}))}{\dim(V(t_1))(t_1-t_1^{-1})}
		\end{align}
		is an invariant of $\mathcal{L}$ which does not depend on the choice of tangle representative.
	\end{thm}
	We refer to this as the $\sltwo$ invariant of $\mathcal{L}$. Also see \cite[Section 6]{Viro} for an alternate proof on the independence of ``cut point.''
\section{Fork Diagrams}\label{sec:quantumgroups}

Here we describe a diagrammatic calculus for the category of representations $V(t)$ which incorporate the inclusions of $V(ts)$ and $V(-ts)$ into $V(t)\otimes V(s)$ described in Theorem \ref{thm:decomp}. These inclusions and their duals are given diagrammatically by signed trivalent vertices which we call elementary forks. We study how forks interact with each other and with the elementary tangle morphisms in Figure \ref{fig:RT}. Recall Conventions \ref{conv:R} and \ref{conv:nondeg}, positive crossings are associated with $R_{ts}$ and all parameters are sufficiently generic.
 
\begin{defn}
	An \textit{elementary fork} is either of the following intertwiners  or their duals:\begin{align}
	\v_+^{t,s}:V(ts)\rightarrow V(t)\otimes V(s)&& \text{or} &&	\v_-^{t,s}:V(-ts)\rightarrow V(t)\otimes V(s),
	\end{align} completely determined by $\v_+^{t,s}(v_0^{ts})=v_0^t\otimes v_0^s$ and $\v_-^{t,s}(v_1^{-ts})=v_1^t\otimes v_1^s$.\\ 

For $\sigma\in\{+,-\}$, the dual $\left(\v_\sigma^{t,s}\right)^*$ is a surjection denoted $\n_\sigma^{t,s}$. These maps are dual in the sense that $
	\n_{\sigma_2}^{t,s}\circ\v_{\sigma_1}^{t,s}=\delta_{\sigma_1\sigma_2}id_{V(\sigma_1ts)}$ for $\sigma_1,\sigma_2\in\{+,-\}$ and where $\delta$ is the Kronecker delta. A \textit{fork} is given by compositions of elementary forks tensored with identity maps. 
\end{defn}

For convenience of the reader, all nonzero values of elementary forks are given below.

\begin{lem}
	Elementary forks are given explicitly by the following evaluations:
	\begin{align}
		&\v_+^{t,s}(v_0^{ts})=v_0^t\otimes v_0^s,
		&&\v_-^{t,s}(v_1^{-ts})=v_1^t\otimes v_1^s,\label{eqn:def}\\
		&\v_+^{t,s}(v_1^{ts})=\Delta(F)(v_0^t\otimes v_0^s), &
		&\v_-^{t,s}(v_0^{-ts})=-\dfrac{1}{\floor{ts}}\Delta(E)(v_1^t\otimes v_1^s),\label{eqn:inc}\\
		&\n_+^{t,s}(v_0^{t}\otimes v_0^s)=v_0^{ts},
		&&\n_-^{t,s}(v_1^{t}\otimes v_1^{s})=v_1^{ts},\label{eqn:ddef}\\
		&\n_+^{t,s}(v_0^t\otimes v_1^s)=\frac{\floor{s}}{\floor{ts}} v_1^{ts}, &	
		&\n_-^{t,s}(v_0^t\otimes v_1^s)= v_0^{-ts},\label{eqn:proj1}\\
		& \n_+^{t,s}(v_1^t\otimes v_0^s)=\frac{s\floor{t}}{\floor{ts}} v_1^{ts}, &
		&\n_-^{t,s}(v_1^t\otimes v_0^s)= -\frac{1}{t}v_0^{-ts}.\label{eqn:proj2}
	\end{align}
\end{lem}
\begin{proof}
	 Recall $\floor{t}=\dfrac{t-t^{-1}}{\zeta-\zeta^{-1}}$. Equations in (\ref{eqn:def}), (\ref{eqn:ddef}) are given and those in (\ref{eqn:inc}) are immediate from the intertwiner property. To evaluate $\n_\sigma^{t,s}(v_0^t\otimes v_1^s)$, write 
	\begin{align*}
	v_0^t\otimes v_1^s=\dfrac{\floor{s}}{\floor{ts}}\Delta(F)(v_0^t\otimes v_0^s)-\dfrac{1}{\floor{s}}\Delta(E)(v_1^t\otimes v_1^s).
	\end{align*}

Evaluating the forks on this expression gives us the equations in (\ref{eqn:proj1}),\begin{align*}
	\n_+^{t,s}(v_0^t\otimes v_1^s)=\n_+^{t,s}\left(\dfrac{\floor{s}}{\floor{ts}}\Delta(F)(v_0^t\otimes v_0^s)\right)=\dfrac{\floor{s}}{\floor{ts}}F\cdot \n_+^{t,s}\left(v_0^t\otimes v_0^s\right)=\dfrac{\floor{s}}{\floor{ts}}v_1^{ts},\\
	\n_-^{t,s}(v_0^t\otimes v_1^s)=\n_-^{t,s}\left(-\dfrac{1}{\floor{ts}}\Delta(E)(v_1^t\otimes v_1^s)\right)=-\dfrac{1}{\floor{ts}}E\cdot\n_-^{t,s}\left(v_1^t\otimes v_1^s\right)=v_0^{-ts}.
\end{align*}
The evaluation of $\n_\sigma^{t,s}$ on $v_1^t\otimes v_0^s=\dfrac{s\floor{t}}{\floor{ts}}\Delta(F)(v_0^t\otimes v_0^s)+\dfrac{1}{t\floor{ts}}\Delta(E)(v_1^t\otimes v_1^s)$ is similar.
\end{proof}

\begin{figure}[h!]
	\includegraphics{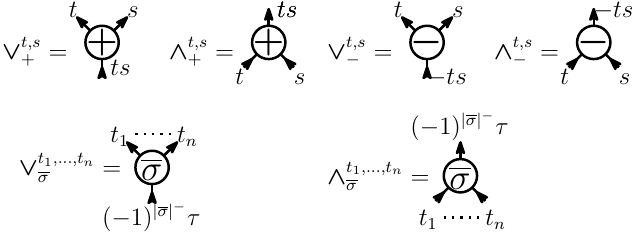}
	\caption{Diagrams for elementary and general forks.}
	\label{fig:elementaryforks}
\end{figure}

We present forks diagrammatically as in Figure  \ref{fig:elementaryforks}. A right-branching fork with $n-1$ vertices is characterized by its sign data $\overline{\sigma}\in \{+,-\}^{n-1}$, which we label from bottom to top. Each such fork determines a unique inclusion ${\v_{\overline{\sigma}}^{t_1,\dots, t_{n}}:V\left(\prod\left(\overline{\sigma}_it_i\right)\right)\to \bigotimes_{i=1}^{n} V(t_i)}$ for some choice of $t_1,\dots, t_{n}$ with dual denoted $\n_{\overline{\sigma}}^{t_1,\dots, t_{n}}$. We use $\tau$ to denote the product $t_1\cdots t_n$. Let $|\overline{\sigma}|^+$ and $|\overline{\sigma}|^-$ denote the total numbers of positive and negative terms in $\overline{\sigma}$, respectively. The sum of all entries in $\overline {\sigma}$ is given by  $\|\overline{\sigma}\|=|\overline{\sigma}|^+-|\overline{\sigma}|^-$ and satisfies $|\overline{\sigma}|^\pm=\dfrac{n-1\pm\|\overline{\sigma}\|}{2}.$\\

We will make use of the ``\textit{channeling identity}''  in the latter parts of this paper to introduce forks into diagrams where they are not already present:
\begin{align}\label{eqn:channeling}
	id_{\bigotimes V(t_i)}=\sum_{\overline{\sigma}} \v_{\overline{\sigma}}^{t_1,\dots, t_n}\circ \n_{\overline{\sigma}}^{t_1,\dots, t_n}.
\end{align}
 The identity is a consequence of the composition $\v_{\overline{\sigma}}^{t_1,\dots, t_n}\circ \n_{\overline{\sigma}}^{t_1,\dots, t_n}$ being a projection onto the direct summand $V\left(\tau(-1)^{|\overline{\sigma}|^-}\right)$ as a submodule of $\bigotimes_{i=1}^{n} V(t_i)$.
 
 \begin{lem}
 	Fix $t_1,\dots, t_n$. The sets of right-branching forks 
 	\begin{align*}
 		\{\v_{\overline{\sigma}}^{t_1,\dots, t_{n}}:(-1)^{|\sigma|^-}=1\} && \text{and} && \{\v_{\overline{\sigma}}^{t_1,\dots, t_{n}}:(-1)^{|\sigma|^-}=-1\}
 	\end{align*}
 are bases of $\Hom(V(\tau),  \bigotimes_{i=1}^{n} V(t_i))$ and $\Hom(V(-\tau),  \bigotimes_{i=1}^{n} V(t_i))$, respectively.
 \end{lem}
\begin{proof}
	Having assumed genericity, according to Theorem \ref{thm:decomp}, 
 \begin{align*}
  \bigotimes_{i=1}^{n} V(t_i)\cong \left(\left(\C(t_1,\dots, t_n)^2\right)^{n-2}\otimes V\left(\tau\right)\right)\oplus \left(\left(\C(t_1,\dots, t_n)^2\right)^{n-2}\otimes V\left(-\tau\right) \right) 
 \end{align*}
 for $n\geq2$. Tensoring with $\left(\C(t_1,\dots, t_n)^2\right)^{n-2}$ here denotes a multiplicity of $2^{n-2}$. By Schur's Lemma, $\Hom(V(\tau),  \bigotimes_{i=1}^{n} V(t_i))$ and $\Hom(V(-\tau),  \bigotimes_{i=1}^{n} V(t_i))$ are vector spaces of dimension $2^{n-2}.$ Varying $\overline{\sigma}$, there are $2^{n-1}$ forks $\v_{\overline{\sigma}}^{t_1,\dots, t_{n}}:V((-1)^{|\sigma|^-}\tau)\to \bigotimes_{i=1}^{n} V(t_i)$. Thus, $2^{n-2}$ of them have $(-1)^{|\sigma|^-}$ of a given sign. Since each of $\{\v_{\overline{\sigma}}^{t_1,\dots, t_{n}}:(-1)^{|\sigma|^-}=1\}$ and $\{\v_{\overline{\sigma}}^{t_1,\dots, t_{n}}:(-1)^{|\sigma|^-}=-1\}$ are linearly independent sets, they are bases of the homomorphism spaces.
\end{proof}

 Although a fork may belong to $\Hom(V(\pm\tau),  \bigotimes_{i=1}^{n} V(t_i))$, it is not necessarily right-branching and may not be immediately comparable to $\v_{\overline{\sigma}}^{t_1,\dots, t_{n}}$. A left-branching fork can be expressed as combination of right-branching forks via a sequence of local ``associator moves.'' 
 
 \begin{lem}
 	The composition of forks satisfies the associator moves given in Figure \ref{fig:assoc}. 
 	\begin{figure}[h!]
 		\includegraphics{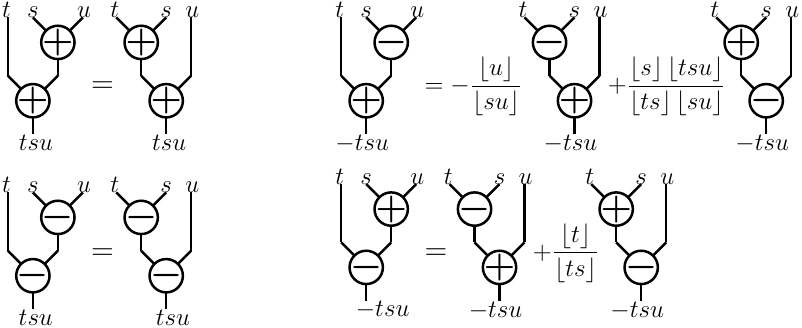}
 		\caption{Associators for forks.}
 		\label{fig:assoc}
 	\end{figure}
 \end{lem}

\begin{proof}
	Recall $\floor{t}=\dfrac{t-t^{-1}}{\zeta-\zeta^{-1}}$. Since forks are intertwiners, it is enough to evaluate all forks in the figure on a single vector. The first left-branching fork is the composition $(\v_+^{t,s}\otimes id_{V(u)})\circ\v_+^{ts,u}$, which evaluates to $v_0^t\otimes v_0^s\otimes v_0^u$ on $v_0^{tsu}$. The right-branching fork $\v_{+,+}^{t,s,u}$ has the same evaluation on $v_0^{tsu}$. Therefore, these two forks are equal. The same argument applies to the equality between $(\v_-^{t,s}\otimes id_{V(u)})\circ\v_-^{-ts,u}$ and $\v_{-,-}^{t,s,u}$.\\
	
	The top-right equality expresses $(\v_-^{t,s}\otimes id_{V(u)})\circ\v_+^{-ts,u}$ as a linear combination of $\v_{+,-}^{t,s,u}$ and $\v_{-,+}^{t,s,u}$. The left-branching fork evaluated on $v_0^{-tsu}$ equals \begin{align*}
		\dfrac{-1}{\floor{ts}}\Delta(E)(v_1^t\otimes v_1^s)\otimes v_0^u=\dfrac{\floor{t}s}{\floor{ts}}v_0^t\otimes v_1^s\otimes v_0^u-\dfrac{\floor{s}}{\floor{ts}}v_1^t\otimes v_0^s\otimes v_0^u.
	\end{align*}
Evaluating the forks on the right side of the equation on $v_0^{-tsu}$ give
\begin{align*}
	\v_{+,-}^{t,s,u}(v_0^{-tsu})&=\dfrac{-1}{\floor{su}}v_0^t\otimes \Delta(E)(v_1^s\otimes v_1^u)=\dfrac{u\floor{s}}{\floor{su}}v_0^t\otimes v_0^s\otimes v_1^u-\dfrac{\floor{u}}{\floor{su}}v_0^t\otimes v_1^s\otimes v_0^u,\\
	\v_{-,+}^{t,s,u}(v_0^{-tsu})&=\dfrac{-1}{\floor{tsu}}(id_{V(t)}\otimes \v_+^{s,u})(\Delta(E)(v_1^{t}\otimes v_1^{su}))
	=(id_{V(t)}\otimes \v_+^{s,u})\left(\dfrac{\floor{t}su}{\floor{tsu}}v_0^{t}\otimes v_1^{su}-\dfrac{\floor{su}}{\floor{tsu}}v_1^{t}\otimes v_0^{su}\right)\\
	&=\dfrac{\floor{t}su}{\floor{tsu}}v_0^{t}\otimes \Delta(F)(v_0^{s}\otimes v_0^u)-\dfrac{\floor{su}}{\floor{tsu}}v_1^{t}\otimes v_0^{s}\otimes v_0^u\\
	&=\dfrac{\floor{t}su}{\floor{tsu}}v_0^{t}\otimes v_1^{s}\otimes v_0^u+\dfrac{\floor{t}u}{\floor{tsu}}v_0^{t}\otimes v_0^{s}\otimes v_1^u-\dfrac{\floor{su}}{\floor{tsu}}v_1^{t}\otimes v_0^{s}\otimes v_0^u.
\end{align*}
Putting these together, it is easy to see
\begin{align*}
	-\dfrac{\floor{t}}{\floor{ts}}\v_{+,-}^{t,s,u}(v_0^{-tsu})
	+\dfrac{\floor{s}\floor{tsu}}{\floor{ts}\floor{su}}\v_{-,+}^{t,s,u}(v_0^{-tsu})
	=(\v_-^{t,s}\otimes id_{V(u)})\circ\v_+^{-ts,u}(v_0^{-tsu}),
\end{align*}
proving the top-right equality.
A similar computation proves the bottom-right equality.
\end{proof}

Computations involving forks and an $R$-matrix can be reduced to associator moves and the diagrammatic manipulations of Figures \ref{fig:forkRmatrix} and \ref{fig:slide}. 

\begin{lem}
	Stacking a diagrammatic crossing on a fork simplifies according to Figure \ref{fig:forkRmatrix}.
\begin{figure}[h!]\centering
	\includegraphics{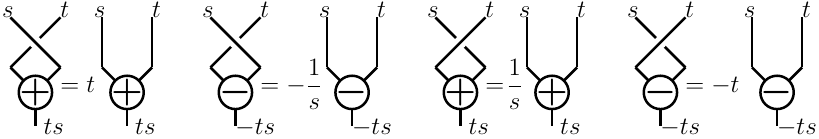}
	\captionof{figure}{Action of $R_{ts}$ on simple forks.}
	\label{fig:forkRmatrix}
\end{figure}
\end{lem}
\begin{proof}
	Recall that diagrammatic crossings denote the normalized $R$-matrix, or its inverse, according to Convention \ref{conv:R}. They are also intertwiners. The lemma follows immediately by evaluating positive signed forks on $v_0^{ts}$ and negative signed forks on $v_1^{-ts}$ using (\ref{eqn:Rts}). 
\end{proof}
\begin{lem}\label{lem:forkcrossingslide}
	Elementary forks slide through crossings according to Figure \ref{fig:slide}.
	\begin{figure}[h!]
		\includegraphics{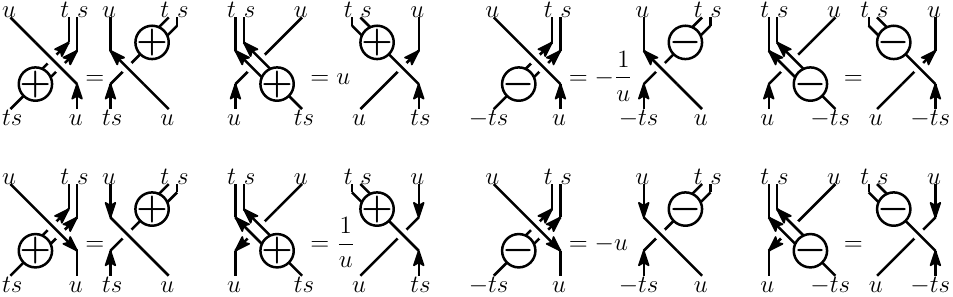}
		\caption{Sliding elementary forks through crossings.}
		\label{fig:slide}
	\end{figure}
\end{lem}
\begin{proof}
	We prove the second equality of the top row by comparing the values of $(id_{V(t)}\otimes R_{us})(R_{ut}\otimes id_{V(s)})(id_{V(u)}\otimes \v_+^{t,s})$ and $(\v_+^{t,s}\otimes id_{V(u)})R_{u,ts}$ on the standard basis vectors of $V(u)\otimes V(ts)$. The other cases are verified by similar arguments. We compute
	 \begin{align*}
		(id_{V(t)}\otimes R_{us})(R_{ut}\otimes id_{V(s)})(id_{V(u)}\otimes \v_+^{t,s})(v_0^u\otimes v_0^{ts})&=
		(id_{V(t)}\otimes R_{us})(R_{tu}\otimes id_{V(s)})(v_0^u\otimes v_0^{t}\otimes v_0^s)
		\\&=u^2v_0^t\otimes v_0^s\otimes v_0^u,
		\\
		(\v_+^{t,s}\otimes id_{V(u)})R_{u,ts}(v_0^u\otimes v_0^{ts})&=(\v_+^{t,s}\otimes id_{V(u)})(uv_0^u\otimes v_0^{ts})\\
		&=u v_0^t\otimes v_0^s\otimes v_0^u,
	\end{align*}
	\begin{align*}
		(id_{V(t)}\otimes R_{us})(R_{ut}\otimes id_{V(s)})&(id_{V(u)}\otimes \v_+^{t,s})(v_0^u\otimes v_1^{ts})
		\\&=(id_{V(t)}\otimes R_{us})(R_{ut}\otimes id_{V(s)})(t^{-1}v_0^u\otimes v_0^t\otimes v_1^s+v_0^u\otimes v_1^t\otimes v_0^s)
		\\&=\dfrac{u}{t}v_0^t\otimes v_1^s\otimes v_0^u+uv_1^t\otimes v_0^s\otimes v_0^u,
		\\
		(\v_+^{t,s}\otimes id_{V(u)})R_{u,ts}(v_0^{u}\otimes v_1^{ts})&=(\v_+^{t,s}\otimes id_{V(u)})(v_1^{ts}\otimes v_0^{u})
		\\&=\dfrac{1}{t}v_0^t\otimes v_1^s\otimes v_0^u+v_1^t\otimes v_0^s\otimes v_0^u,
	\end{align*}
	\begin{align*}
		(id_{V(t)}\otimes R_{us})(R_{ut}\otimes id_{V(s)})(id_{V(u)}\otimes &\v_+^{t,s})(v_1^u\otimes v_0^{ts})
		=(id_{V(t)}\otimes R_{us})(R_{ut}\otimes id_{V(s)})(v_1^u\otimes v_0^t\otimes v_0^s)
		\\&=\dfrac{u^2}{ts}v_0^u\otimes v_0^t\otimes v_1^s+\dfrac{u^2-1}{t}v_0^u\otimes v_1^t\otimes v_0^s+(u^2-1)v_1^u\otimes v_0^t\otimes v_0^s,
		\\
		(\v_+^{t,s}\otimes id_{V(u)})R_{u,ts}(v_0^u\otimes v_1^{ts})&=(\v_+^{t,s}\otimes id_{V(u)})(v_1^{ts}\otimes v_0^{u})
		\\&=\dfrac{u}{ts}v_0^t\otimes v_0^s\otimes v_1^u+\dfrac{u^2-1}{ut}v_0^t\otimes v_1^s\otimes v_0^u+\dfrac{(u^2-1)}{u}v_1^t\otimes v_0^s\otimes v_0^u,
	\end{align*}
	\begin{align*}\allowdisplaybreaks
		(id_{V(t)}\otimes R_{us})(R_{ut}\otimes id_{V(s)})(id_{V(u)}\otimes &\v_+^{t,s})(v_1^u\otimes v_1^{ts})
		\\&=(id_{V(t)}\otimes R_{us})(R_{ut}\otimes id_{V(s)})(t^{-1}v_1^u\otimes v_0^t\otimes v_1^s+v_1^u\otimes v_1^t\otimes v_0^s)
		\\&=\dfrac{-u}{t^2s}v_0^t\otimes v_1^s\otimes v_1^u-\dfrac{u}{ts}v_1^t\otimes v_0^s\otimes v_1^u,
		\\
		(\v_+^{t,s}\otimes id_{V(u)})R_{u,ts}(v_1^u\otimes v_1^{ts})&=(\v_+^{t,s}\otimes id_{V(u)})(-(ts)^{-1}v_1^{ts}\otimes v_1^{u})
		\\&=\dfrac{-1}{t^2s}v_0^t\otimes v_1^s\otimes v_1^u-\dfrac{1}{ts}v_1^t\otimes v_0^s\otimes v_1^u.
	\end{align*}
This shows that $(id_{V(t)}\otimes R_{us})(R_{ut}\otimes id_{V(s)})(id_{V(u)}\otimes \v_+^{t,s})=u(\v_+^{t,s}\otimes id_{V(u)})R_{u,ts}$.
\end{proof}
{\begin{rem}\label{rem:notnatural}
		Passing the fork through a \textit{normalized} crossing is not natural, which is why there are factors appearing in Lemma \ref{lem:forkcrossingslide}.
\end{rem}}

\begin{lem}\label{lem:slide}
	 For any inclusion of $V\left((-1)^{|\overline{\sigma}|^-} \tau\right)$ into $\bigotimes_{i=1}^n V(t_i)$ by a right-branching fork diagram with vertex labels $\overline{\sigma}_i\in\{\pm1\}$, we have the  equalities of Figure \ref{fig:bigslide}. The product of the factors $\left(\dfrac{(-1)^{n}}{\tau}\right)^{{|\overline{\sigma}|}^-}$ and $\left((-1)^{|\overline{\sigma}|^-}\tau\right)^{|\overline{\sigma}|^+}$ that appear equals $\tau^{\|\overline{\sigma}\|}$.

	 \allowdisplaybreaks
	 \begin{figure}[h!]
	 	\centering
	 	\includegraphics{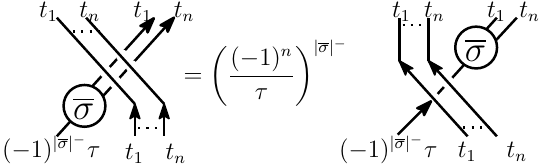}
	 	
	 	\vspace{2em}
	 		
	 \includegraphics{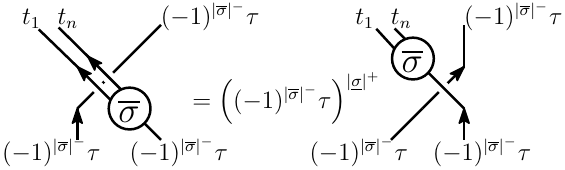}\caption{Equalities of Lemma \ref{lem:slide}.}\label{fig:bigslide}
 \end{figure}  
\end{lem}

\begin{proof}
	By repeated application of the third diagram in Figure \ref{fig:slide}, passing a negative signed fork under all the strands yields a factor $(-1)^n\tau^{-1}$. Positive signed forks slide freely under the crossing. Thus,  sliding the $\overline{\sigma}$ labeled fork through yields $(-1)^{n|\overline{\sigma}|^-}\tau^{-|\overline{\sigma}|^-}$, $|\overline{\sigma}|^-$ denoting the total number of negative signs appearing in $\overline{\sigma}$.\\
	
	Each positive elementary fork passing over the strand labeled by $(-1)^{|\overline{\sigma}|^-} \tau$  contributes a factor equal to the labeling. Sliding the entire fork over yields $((-1)^{|\overline{\sigma}|^-} \tau)^{|\overline{\sigma}|^+}$. The product of the two terms is
	\begin{align*}
		\left((-1)^{n |\overline{\sigma}|^-}\tau^{-|\overline{\sigma}|^-}\right)\left((-1)^{|\overline{\sigma}|^-\cdot|\overline{\sigma}|^+}\tau^{|\overline{\sigma}|^+}\right)
		&=(-1)^{|\overline{\sigma}|^-(n+|\overline{\sigma}|^+)}\tau^{|\overline{\sigma}|^+-|\overline{\sigma}|^-}\\
		&=(-1)^{|\overline{\sigma}|^-(n+(n-1-|\overline{\sigma}|^-))}\tau^{\|\overline{\sigma}\|}
		=\tau^{\|\overline{\sigma}\|}.\qedhere
	\end{align*}
\end{proof}

\begin{figure}[h!]
  	\includegraphics[scale=1]{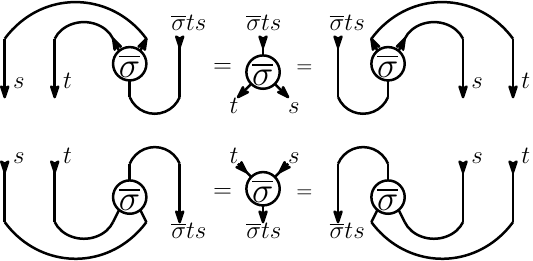}
  	\caption{Downward oriented forks defined using evaluation/coevaluation.}\label{fig:rotated}
  \end{figure}

\begin{figure}[h!]
  	\includegraphics[scale=1]{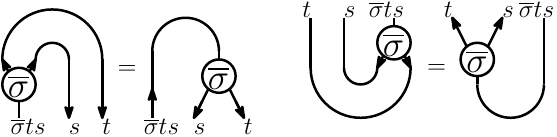}
  	\caption{Forks slide through extrema.}
  	\label{fig:maxmin}
  \end{figure}

Using the evaluation and coevaluation maps we construct downward oriented forks, see Figure \ref{fig:rotated}. Together with equation (\ref{eq:zigzag1}) and (\ref{eq:zigzag2}), Figure \ref{fig:rotated} implies the relations of Figure \ref{fig:maxmin}.
  

\section{Recovering the Conway Potential Function}\label{sec:CPFproof}
\begin{figure}[h!]
	\[\nabla\begin{pmatrix}
		\includegraphics[scale=.9]{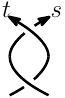} 
	\end{pmatrix}+\nabla\begin{pmatrix}
		\includegraphics[scale=.9]{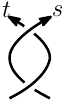} 
	\end{pmatrix}=(ts+t^{-1}s^{-1})\nabla\begin{pmatrix}
		\includegraphics[scale=.9]{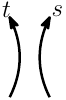} 
	\end{pmatrix}\label{eq:II}\tag{II}\]
	\begin{align*}
		(t^{-1}s^{-1}-ts)&\left\{\nabla\begin{pmatrix}
			\includegraphics[scale=1]{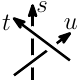} 
		\end{pmatrix}+\nabla\begin{pmatrix}
			\includegraphics[scale=1]{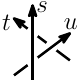} 
		\end{pmatrix}\right\}\\
		+(s^{-1}u^{-1}-su)&\left\{\nabla\begin{pmatrix}
			\includegraphics[scale=1]{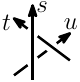} 
		\end{pmatrix}+\nabla\begin{pmatrix}
			\includegraphics[scale=1]{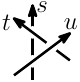} 
		\end{pmatrix}\right\}
		\\
		+(tu^{-1}-t^{-1}u)&\left\{\nabla\begin{pmatrix}
			\includegraphics[scale=1]{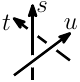}\end{pmatrix}+\nabla\begin{pmatrix}
			\includegraphics[scale=1]{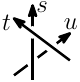} 
		\end{pmatrix}\right\}=0\label{eq:III}\tag{III}
	\end{align*}
	\[\nabla\begin{pmatrix}
		~~\includegraphics[scale=.9]{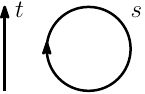} 
	\end{pmatrix}=0\label{eq:IO}\tag{IO}\]
	\[
	\nabla\begin{pmatrix}
		\includegraphics[scale=.9]{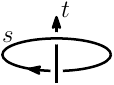} 
	\end{pmatrix}=(t-t^{-1})\nabla\begin{pmatrix}
		~~\includegraphics[scale=.9]{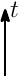} 
	\end{pmatrix}\label{eq:Phi}\tag{$\Phi$}\]
	\[
	\nabla\begin{pmatrix}
		\includegraphics[scale=.9]{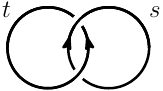} 
	\end{pmatrix}=1\label{eq:H}\tag{H}
	\]\caption{Jiang's skein relations for the Conway Potential Function.}\label{fig:rels}
\end{figure}

In this section, we prove the $\sltwo$ invariant satisfies Jiang's axioms for the CPF in Figure \ref{fig:rels}. Lemma \ref{lem:factor}, stated at the end of this section, shows how signs factor out of the CPF and is essential in the proofs of Theorems \ref{thm:sat} and \ref{thm:link}.

\begin{proof}[Proof of  \emph{(\ref{eq:II})}]
	 Jiang characterizes both (\ref{eq:II}) and (\ref{eq:III}) in the context of braids, and so it is enough to prove the equalities only for the diagrams shown. Relation (\ref{eq:II}) translates to
	\begin{align*}
	R_{st}R_{ts}+(R_{st}R_{ts})^{-1}=\left(ts+\dfrac{1}{ts}\right)\cdot id_{V(t)}\otimes id_{V(s)}.
	\end{align*} Then as 
	\begin{align*}
	R_{st}R_{ts}=\arraycolsep=1pt
	\begin{bmatrix}
	st&0&0&0\\ \noalign{\medskip}0&{\dfrac {s
		}{t}}&s-{\dfrac {s}{{t}^{2}}}&0\\ \noalign{\medskip}0&s-\dfrac{1}{s}&st-{
		\dfrac {s}{t}}+{\dfrac {1}{st}}&0\\ \noalign{\medskip}0&0&0&{\dfrac {1}{s
			t}}
	\end{bmatrix} &&\text{and} && (R_{st}R_{ts})^{-1}=
	\begin{bmatrix}\arraycolsep=1pt
	{\dfrac {1}{st}}&0&0 & 0
	\\ \noalign{\medskip}0&st-{\dfrac {s}{t}}+{\dfrac {1}{st}}&{\dfrac {s
		}{{t}^{2}}}-s&0\\ \noalign{\medskip}0&\dfrac{1}{s}-s&{\dfrac {s}{t}}&0
	\\ \noalign{\medskip}0&0&0&st
	\end{bmatrix},\end{align*} 
	the relation is seen to hold. 
\end{proof}
\begin{proof}[Proof of \emph{(\ref{eq:III})}]
	Relation (\ref{eq:III}) is the most complicated relation to check, which is a relation between six maps $:V(u)\otimes V(s)\otimes V(t)\rightarrow V(t)\otimes V(s)\otimes V(u).$	We consider an inclusion of $V(-tsu)$ into $V(t)\otimes V(s)\otimes V(u)$. We demonstrate (\ref{eq:III}) holds on this summand and a similar verification can be made for the other inclusion of $V(-tsu)$ and those of $V(tsu)$. We resolve the diagrams into simple forks using the associator and $R$-matrix moves of Figures \ref{fig:assoc} and \ref{fig:forkRmatrix}. At the $i$-th diagrammatic crossing define $C_i^\pm(\sigma,t,s)$ to be the coefficient obtained in Figure \ref{fig:forkRmatrix} from the designated labels, with $\pm$ indicating  the sign of the crossing. We compute the action of a left-branching fork on the braids in (\ref{eq:III}) using Figure \ref{fig:proofIII}.
	
	\begin{figure}[h!]
		\includegraphics{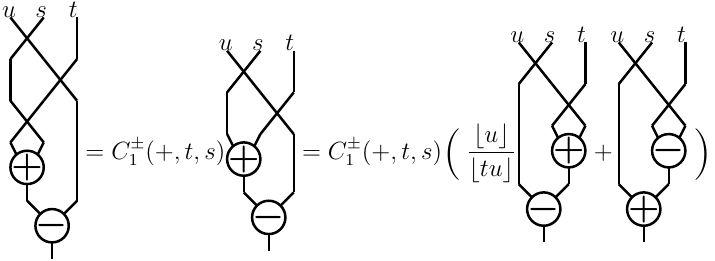}\\

		\includegraphics{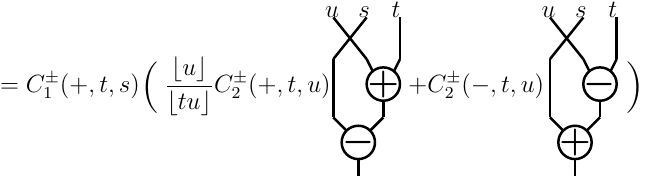}\\
		
		\includegraphics{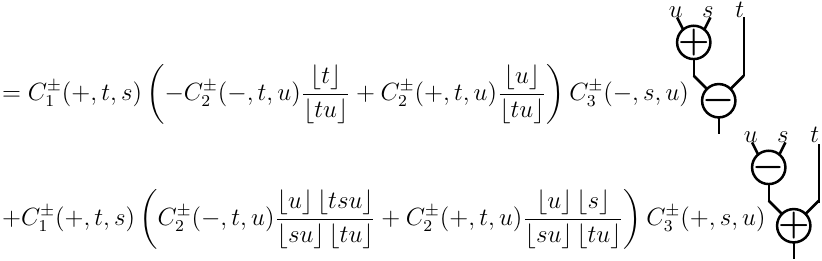}
		\caption{Applying the $(+,-)$ fork to the braids for relation (\ref{eq:III}).}\label{fig:proofIII}
	\end{figure}

	 For the resulting left-branching fork $\overline{\sigma}=(-,+)$, we sum over the different crossings in (\ref{eq:III}),
	
\begin{gather*}
	\centering
		\begin{aligned}
		-&(t^{-1}s^{-1}-ts)\left(\frac{\floor{u}}{\floor{tu}}\left(\frac{t^2}{s}+\frac{1}{u^2s}\right)+\frac{\floor{t}}{\floor{tu}}\left(\frac{ts}{u}+\frac{t}{su}\right)\right)\\
		-&(s^{-1}u^{-1}-us)\left(\frac{\floor{u}}{\floor{tu}}\left(\frac{t}{su}+\frac{ts}{u}\right)+\frac{\floor{t}}{\floor{tu}}\left(\frac{1}{u^2s}+t^2s\right)\right)\\
		-&(tu^{-1}-ut^{-1})\left(\frac{\floor{u}}{\floor{tu}}\left(\frac{t^2}{u}+\frac{1}{u}\right)+\frac{\floor{t}}{\floor{tu}}\left(\frac{t}{u^2}+t\right)\right)=0.
	\end{aligned}
\end{gather*}

Similarly for the $\overline{\sigma}=(+,-)$ fork, after factoring out $\dfrac{\floor{u}}{\floor{su}\floor{tu}}$, 
\begin{align*}
	&(t^{-1}s^{-1}-ts)\left(-\floor{tsu}\left(\frac{t}{u^2}+t\right)+\floor{s}\left(\frac{t^2}{u}+\frac{1}{u}\right)\right)\\
	+&(s^{-1}u^{-1}-us)\left(-\floor{tsu}\left(\frac{t^2}{u}+\frac{1}{u}\right)+\floor{s}\left(\frac{t}{u^2}+t\right)\right)\\
	+&(tu^{-1}-ut^{-1})\left(-\floor{tsu}\left(\frac{ts}{u}+\frac{t}{su}\right)+\floor{s}\left(t^2s+\frac{1}{su^2}\right)\right)=0.\qedhere
\end{align*}
\end{proof}

\begin{proof}[Proof of \emph{(\ref{eq:IO})}]	The map defined by the circle colored by $t$ is multiplication by the quantum dimension of $V(t)$. Since
this dimension is zero,
the invariant assigns zero to such links.
\end{proof} 

\begin{proof}[Proof of \emph{(\ref{eq:Phi})}]
	As a $V(t)$ endomorphism, equation $(\Phi)$ determines the map \begin{align*}
	&(id_{V(t)}\otimes \widetilde{ev}_{V(s)})\cdot 
	((R_{st}R_{ts})^{-1}\otimes id_{V(s)})\cdot
	(id_{V(t)}\otimes coev_{V(s)})=(t-t^{-1})id_{V(t)}. \hfill\qedhere
	\end{align*}
\end{proof}

\begin{proof}[Proof of \emph{(\ref{eq:H})}]
	Cutting either string in the Hopf link produces the 1-tangle seen in (\ref{eq:Phi}). Thus, the $V(t)$ endomorphism associated to the Hopf link having cut the $i$-colored component is ${(t-t^{-1})id_{V(t)}}$. Tracing off this strand and dividing by $\dim{V(t)}({t-t^{-1}})$ yields 1.
\end{proof}

We refer to Lemma 4.1(c) of \cite{Melikhov} to measure changes in the CPF at negative parameters.

\begin{lem}\label{lem:factor}
	Let $\mathcal{L}$ be a link with $n\geq 1$ components and linking numbers $l_{jk}$. Suppose the number of components colored by $t_i$ is $m_i$. Then
	\begin{align*}
		\nabla(\mathcal{L})(t_1,\dots,-t_i,\dots, t_n)=(-1)^{m_i+p}\nabla(\mathcal{L})(t_1,\dots, t_n),
	\end{align*}
with $p$ equal to the sum over all $l_{jk}$ for which $t_{c_j}=t_i$ and $t_{c_k}\neq t_i$.
\end{lem}
\begin{proof}
	Since the Jiang relations completely determine the CPF, if each relation is compatible with the proposed sign factorization, then it is true for the CPF itself. It is necessary to consider sign changes on any component and repeated colorings.\newline
	
	The relations are readily verified, we discuss one case of relation (\ref{eq:III}). If $u$ is replaced by $-u$, we normalize the number of $u$ colored components and linkings according to the first diagram in (\ref{eq:III}). If $t$ and $s$ are distinct from $u$, the linking between the $u$ stand and the others both change by one in the second diagram. In the remaining diagrams, the linking only changes between either $t$ or $s$ and an extra sign comes from the coefficient that depends on $u$. Therefore, the relation holds in this case. 
\end{proof}
\section{Seifert-Torres Type Formulas}
In this section we prove our main results, Theorems \ref{thm:sat} and \ref{thm:link}. A special case of Theorem \ref{thm:link} for knots is given in Corollary \ref{cor:cable}. \newline

\begin{figure}[h!]
	\centering
		\includegraphics{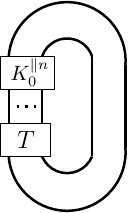}
	\caption{A diagram for a {$(\mathcal{K},\mathcal{T})$-satellite} link.} \label{fig:KTsat}
\end{figure}

\begin{defn}
	Let $\mathcal{K}$ be a 0-framed knot and $\mathcal{T}$ a framed $n$-tangle. We denote by $\mathcal{K}_0$ any framed $1$-tangle with closure equal to $\mathcal{K}$, and $\mathcal{K}_0^{\parallel n}$ to be the $n$-fold parallel cabling of $\mathcal{K}_0$ with respect to its framing. Then the $(\mathcal{K},\mathcal{T})$-\textit{satellite} is given by Figure \ref{fig:KTsat}.
\end{defn} 

The same construction applies to framed links without self-linking on each component. The $(\mathcal{L},\mathcal{T}_1,\dots, \mathcal{T}_n)$-\textit{satellite link} is obtained by forming an $(\mathcal{L}_i,\mathcal{T}_i)$-satellite on each component of $\mathcal{L}$, given tangles $\mathcal{T}_i$.\\

A $(\mathcal{L},\mathcal{T}_1,\dots, \mathcal{T}_n)$-satellite link is an example of a splice link in the sense of \cite[Chapter 1]{EN}. Note that there is a natural inclusion of $\widehat{\mathcal{T}}_i$ into a torus such that the linking number between the meridian $\mu$ and $\widehat{\mathcal{T}}_i$ equals the number of components of $\mathcal{T}_i$. Let $\mathcal{S}_i$ be the link $\widehat{\mathcal{T}_i}\cup S^1$, where $S^1$ is an unknot homologous to $\mu$. Then $\mathcal{L}'$ is obtained by splicing each component of $\mathcal{L}$ with the unknot in the corresponding ${\mathcal{S}_i}$.
	
	\begin{rem}
		The 0-framed condition on $\mathcal{K}$ in the definition of the $(\mathcal{K},\mathcal{T})$-satellite produces a well-defined diagram in terms of abstract knots and simplifies computations. Writhe can be introduced into the diagram by adding full twists to $\mathcal{T}$. 
	\end{rem}

\begin{figure}[h!]
			\includegraphics{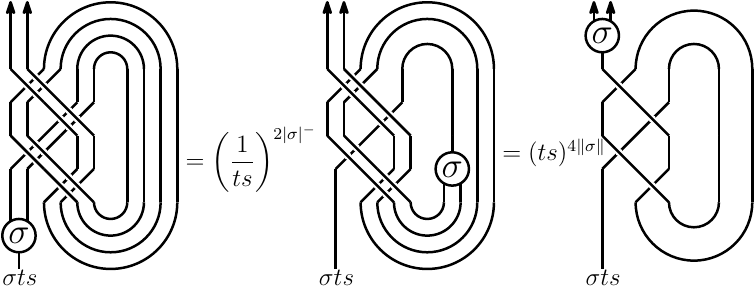}
			\caption{Sliding a fork through an open 2-cabled trefoil.}\label{fig:slideex}
		\end{figure}

	\begin{lem}\label{lem:DK}
		Let $\mathcal{K}_0$ be a 1-component 0-framed 1-tangle obtained from a braid closure. Fix sign data $\overline{\sigma}\in\{\pm1\}^{n-1}$ and set $\tau=t_1\cdots t_n$. The  diagram in Figure \ref{fig:DK} commutes.\vskip.4\baselineskip
		\begin{figure}[h!]
			\centering
			\begin{tikzcd}
			V((-1)^{|\overline{\sigma}|^-} \tau ) \arrow[r, "\bracks{\mathcal{K}_0}"] \arrow[d,"\v_{\overline{\sigma}}^{t_1,\dots, t_n}",swap,hook]
			& 	V((-1)^{|\overline{\sigma}|^-} \tau) \arrow[d, "\v_{\overline{\sigma}}^{t_1,\dots, t_n}" ,hook] \\
			\bigotimes_{i=1}^{n} V(t_i) \arrow[r, "\bracks{\mathcal{K}_0^{\parallel n}}"]
			& \bigotimes_{i=1}^{n} V(t_i)
		\end{tikzcd}
		\caption{Commutative diagram of Lemma \ref{lem:DK}}\label{fig:DK}
		\end{figure}
	\end{lem}
\begin{proof}
Recall that  $\bracks{\mathcal{K}_0}$ and $\bracks{\mathcal{K}_0^{\parallel n}}$ are the endomorphisms of $V((-1)^{|\overline{\sigma}|^-} \tau)$ and $\bigotimes_{i=1}^{n} V(t_i)$ assigned to ${\mathcal{K}_0}$ and ${\mathcal{K}_0^{\parallel n}}$ under the Reshetikhin-Turaev functor. 
	To show the equality, we slide the fork $\v_{\overline{\sigma}}^{t_1,\dots, t_n}$ through $\bracks{\mathcal{K}_0^{\parallel n}}$ to produce $\bracks{\mathcal{K}_0}$ colored by $V((-1)^{|\overline{\sigma}|^-} \tau )$ composed with $\v_{\overline{\sigma}}^{t_1,\dots, t_n}$. We verify that sliding the fork through $\mathcal{K}_0^{\parallel n}$ does not introduce additional factors. Equivalently, we prove the relation in Figure \ref{fig:natural}.
	
\begin{figure}[h!]
 	\centering\includegraphics{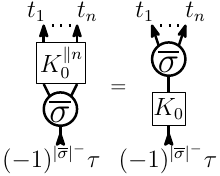}
		\caption{Naturality for writhe-less diagrams.}\label{fig:natural}
 \end{figure}
	
	Since $\mathcal{K}_0^{\parallel n}$ is the cabling of a 1-component tangle, the fork will appear on each side of every cabled crossing as it slides through the diagram. According to Figure \ref{fig:maxmin}, forks slide through extrema. Thus, all cabled strands can be replaced with single strands colored by  $V((-1)^{|\overline{\sigma}|^-} \tau)$ except at crossings where a factor of $\tau^{\pm\|\overline{\sigma}\|}$ is introduced depending on the sign of the crossing, as in Lemma \ref{lem:slide}. Our assumption of $\mathcal{K}_0$ being obtained from a braid presentation implies both strands at each crossing are oriented upwards, and so all factors are of the required form. Since $\mathcal{K}_0$ has zero writhe, the overall contribution of these factors is trivial. Figure \ref{fig:slideex} depicts the intermediate steps of the slide on an open 2-cabled trefoil with \textit{nonzero} writhe.
\end{proof}

\begin{rem}\label{rem:slidelink}
	We have a similar result if the 1-tangle in Lemma \ref{lem:DK} has multiple components. The additional factors which appear are $\prod_j \left((-1)^{|\overline{\sigma}|^-}t_{j}^{\|\overline{\sigma}\|}\right)^{l_{1j}}$, where $l_{1j}$ is the linking between the open component and component $j$. This can be readily computed from Lemma \ref{lem:slide}.
\end{rem}
\begin{lem}\label{lem:h}
	For any sign data $\overline{\sigma}\in\{\pm1\}^{n-1}$ and $t_1,\dots, t_n\in\C^\times$,
	\begin{align*}
		\tr(\n_{\overline{\sigma}}^{t_1,\dots, t_n}\circ(id_{V(t_1)}\otimes \rho_{t_2,\cdots, t_n}(\Delta^{n-2}(K^{-1}))\circ\v_{\overline{\sigma}}^{t_1,\dots, t_n})=
		(-1)^{|\overline{\sigma}|^-}\dfrac{2(t_1-t_1^{-1})}{\tau-\tau^{-1}}.
	\end{align*}
\end{lem}
\begin{proof}
	Diagrammatically, the action of $K^{-1}$ on a tensor factor is denoted by a bead on the associated strand. Since each fork is an intertwiner, $\rho_{t_1,t_2}(\Delta(K^{-1}))\circ\v^{t_1,t_2}_{\sigma}=\v^{t_1,t_2}_{\sigma}\circ\rho_{t_1t_2}(K^{-1})$ for $\sigma\in\{\pm1\}$. This is shown diagrammatically in Figure \ref{fig:Kslide}. 
	
	\begin{figure}[h!]
		\centering 
		\includegraphics{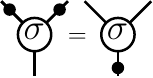}
		\caption{Sliding $\Delta(K^{-1})$ through a fork.}\label{fig:Kslide}
	\end{figure}
	
	Since forks are right-branching, we iterate the above and by duality, we have
	\begin{align*}
		\n_{\overline{\sigma}}^{t_1,\dots, t_n}\circ(id_{V(t_1)}\otimes \rho_{t_2,\cdots, t_n}(\Delta^{n-2}(K^{-1}))\circ\v_{\overline{\sigma}}^{t_1,\dots, t_n}\\
		=\n_{\overline{\sigma}_1}^{t_1,s}\circ(id_{V(t_1)}\otimes \rho_{s}(K^{-1}))\circ\v_{\overline{\sigma}_1}^{t_1,s},
	\end{align*}
	where $s=\overline{\sigma}_1(-1)^{|\overline{\sigma}|^-}t_2\cdots t_n$. For $\overline{\sigma}_1=+1$, $s=(-1)^{|\overline{\sigma}|^-}t_1^{-1}\tau$ and
	\begin{align*}
		&\n_{\overline{\sigma}_1}^{t_1,s}\circ(id_{V(t_1)}\otimes \rho_{s}(K^{-1}))\circ\v_{\overline{\sigma}_1}^{t_1,s}(v_0^{\tau})=\dfrac{1}{s}v_0^\tau,\\
		&\n_{\overline{\sigma}_1}^{t_1,s}\circ(id_{V(t_1)}\otimes \rho_{s}(K^{-1}))\circ\v_{\overline{\sigma}_1}^{t_1,s}(v_1^{\tau})=
		\n_{\overline{\sigma}_1}^{t_1,s}(\dfrac{1}{s}v_1^{t_1}\otimes v_0^s-\dfrac{1}{t_1s}v_0^{t_1}\otimes v_1^s)
		=\dfrac{\floor{t_1}}{\floor{t_1s}}-\dfrac{\floor{s}}{\floor{t_1s}t_1s}v_1^\tau,\\
		&\tr(\n_{\overline{\sigma}_1}^{t_1,s}\circ(id_{V(t_1)}\otimes \rho_{s}(K^{-1}))\circ\v_{\overline{\sigma}_1}^{t_1,s})=(-1)^{|\overline{\sigma}|^-}\left(
		\dfrac{t_1}{\tau}+\dfrac{\floor{t_1}}{\floor{\tau}}-\dfrac{\floor{\tau t_1^{-1}}}{\floor{\tau}\tau}
		\right)=(-1)^{|\overline{\sigma}|^-}\dfrac{2(t_1-t_1^{-1})}{\tau-\tau^{-1}}.
	\end{align*}
The computation for $\overline{\sigma}_1=-1$ is similar.
\end{proof}

\begin{lem}\label{lem:fork}
	Let $\mathcal{T}$ be a string link with coloring $\overline{c}$ valued in $\{t_1,\dots, t_m\}^n$ such that $\widehat{\mathcal{T}}$ has a well defined coloring. Set $\tau={c_1}\cdots {c_n}$. Then \begin{align*}
		\nabla(\widehat{\mathcal{T}})(t_1,\dots, t_m)=\dfrac{1}{2(\tau-\tau^{-1})}\sum_{\overline{\sigma}}(-1)^{|\overline{\sigma}|^-}\tr\left(\n_{\overline{\sigma}}^{{c_1},\dots, {c_n}}\circ \bracks{\mathcal{T}}\circ \v_{\overline{\sigma}}^{{c_1},\dots, {c_n}}\right).
	\end{align*}
\end{lem}
\begin{proof}
	By Theorem \ref{thm}, $\nabla(\widehat{\mathcal{T}})$ is given by $\tr(id_{V(c_1)}\otimes \rho_{c_2,\cdots, c_n}(\Delta^{n-2}(K^{-1}))\circ\bracks{\mathcal{T}})$ normalized by $2(c_{1}-c_1^{-1})$. We sum over pairs of forks labeled by $\overline{\sigma}$ and $\overline{\sigma}'$ directly above and below $\mathcal{T}$ by two applications of the ``channeling identity'' in (\ref{eqn:channeling}). Sliding the bottom fork along the closure of the diagram allows us to write the diagram as 
	\begin{align*}
		\tr\left(\n_{\overline{\sigma}}^{c_1,\dots, c_n}\circ (id_{V(c_1)}\otimes \rho_{c_2,\cdots, c_n}(\Delta^{n-2}(K^{-1}))\circ \v_{\overline{\sigma}'}^{c_1,\dots, c_n}\circ \n_{\overline{\sigma}'}^{c_1,\dots, c_n}\circ \bracks{\mathcal{T}}\circ \v_{\overline{\sigma}}^{c_1,\dots, c_n}\right).
	\end{align*}
	Since $\n_{\overline{\sigma}}^{c_1,\dots, c_n}\circ (id_{V(c_1)}\otimes \rho_{c_2,\cdots, c_n}(\Delta^{n-2}(K^{-1})))\circ \v_{\overline{\sigma}'}^{c_1,\dots, c_n}
		$ equals $\delta_{\overline{\sigma},{\overline{\sigma}'}}(-1)^{|\overline{\sigma}|^-}\dfrac{c_1-c_1^{-1}}{\tau-\tau^{-1}}$ by orthogonality and Lemma \ref{lem:h}, we have the desired result after normalizing.
\end{proof}
	
\sat
{\begin{proof}
	Diagrams for this proof are given in Figure \ref{fig:satproof}, where $n=2$ and forks are labeled by $\overline{\sigma}$ and $\overline{\mu}$. We assume that $\mathcal{L}$ is given by the closure of a braid and $\mathcal{L}_0$ is an $n$-component string link obtained from its partial closure. The cabled string link $\mathcal{L}_0^{\parallel \overline{d}}$ is then stacked above $\mathcal{T}_1,\dots,\mathcal{T}_n$. We apply the ``channeling identity'' above each $\mathcal{T}_i$ and slide each fork through the cabled link diagram. We label the sign data of the $i$-th fork $\overline{\sigma}^i$. Since $\mathcal{L}$ has zero linking matrix, no additional factors appear according to Lemma \ref{lem:slide} and Remark \ref{rem:slidelink}.\\
	
	The diagrams $\mathcal{T}_2,\dots, \mathcal{T}_n$  sandwiched by forks factor out. We also take their canonical trace and divide by $2$. The diagram for $\mathcal{T}_1$ does not factor as the first strand is not acted on by $K^{-1}$. It follows that a 1-tangle representative with closure $\mathcal{L}$ also factors. This 1-tangle acts by the scalar $\nabla(\mathcal{L})\left((-1)^{|\overline{\sigma}^1|^-}\tau_1,\dots, (-1)^{|\overline{\sigma}^n|^-}\tau_n\right)(-1)^{|\overline{\sigma}^1|^-}(\tau_1-\tau_1^{-1})$. By Lemma \ref{lem:factor}, each sign factors from the CPF to yield $\nabla(\mathcal{L})(\tau_1,\dots,\tau_n)(\tau_1-\tau_1^{-1})(-1)^{\sum_{i=2}^n|\overline{\sigma}^i|^-}$.\\
	
	After factoring the diagram for $\mathcal{L}$, $\bracks{\mathcal{T}_1}$ composed with the channeling identity over $\overline{\sigma}^1$ and $id_{V(c_{11})}\otimes\rho_{c_{12},\dots, c_{1d_1}}\left(\Delta^{d_1-1}(K^{-1})\right)$ remains. We may sum over the forks around each of the tangles $\mathcal{T}_1,\dots, \mathcal{T}_n$ independently. Summing over $\overline{\sigma}^1$ yields $2(t_1-t_1^{-1})\nabla(\mathcal{T}_1)$, the coefficient cancels after normalization.  We apply Lemma \ref{lem:fork} for the remaining tangles and signs $(-1)^{|\overline{\sigma}^i|^-}$ and each contributes $2(\tau_i-\tau_i^{-1})\nabla(\widehat{\mathcal{T}_i})(t_{i,1},\dots, t_{i,d_i})$. The factors of $2$ cancel with the earlier terms from the canonical trace. Note that $\nabla(\mathcal{T}_1)$ pairs with $(\tau_1-\tau_1^{-1})$ from the earlier steps involving $\nabla(\mathcal{L})$, which produces the desired formula. \qedhere

\begin{figure}[h!]
	\centering\includegraphics{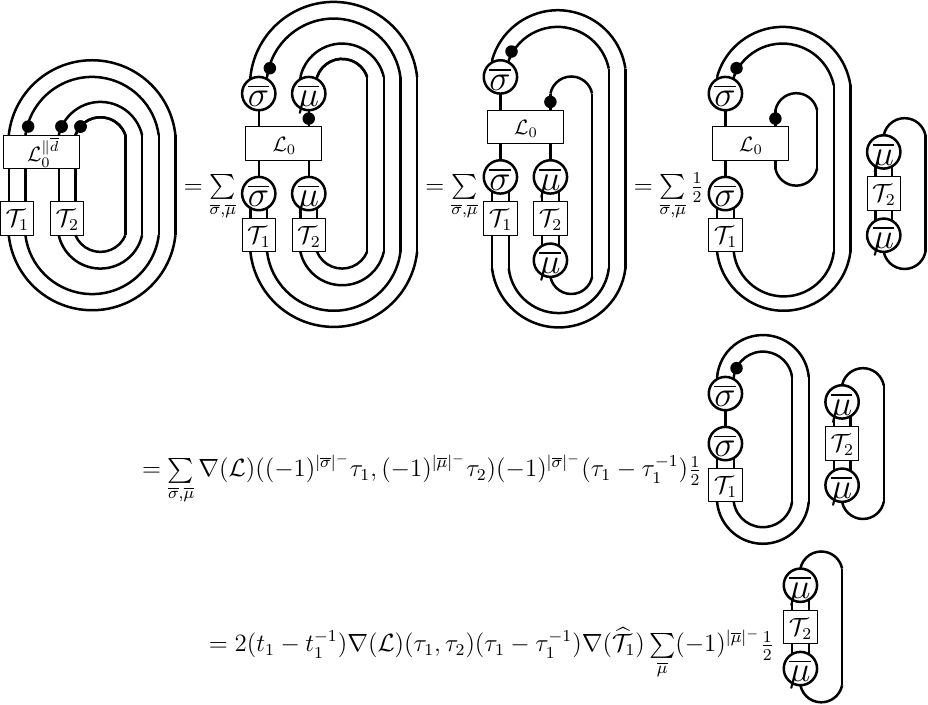}
	\caption{Proof of Theorem \ref{thm:sat} where $\mathcal{L}$ has two components and are forks are labeled by $\overline{\sigma}$ and $\overline{\mu}$.}\label{fig:satproof}
\end{figure}
\end{proof}}

\begin{rem}
	The theorem relies on $\mathcal{T}$ being a string link. Since $P(-1)\cong V(t)\otimes V(t)^*$ is indecomposable, $id_{P(-1)}$ is not expressible as a sum over forks. In which case, our argument is no longer valid.
\end{rem}

\link

\begin{proof}
	We consider presentation of $\mathcal{L}$ as the closure of a braid $\mathcal{L}_0$. On each of the $i$ components of $\mathcal{L}^{\parallel \overline{d}}$, we take a sum over all right-branching forks with sign data $\overline{\sigma}^i$ via the ``channeling identity.'' Sliding all forks through the cabled diagram yields $\widehat{\mathcal{L}_0}$ colored by $(-1)^{|\overline{\sigma}^i|^-}\tau_i$ on its $i$-th component. By Lemma \ref{lem:slide} and Remark \ref{rem:slidelink}, this sliding move also produces an overall factor $\tau_i^{\|\overline{\sigma}^i\|}$ on the $i$-th component from writhe, and $(-1)^{|\overline{\sigma}^i|^-}\tau_j^{\|\overline{\sigma}^i\|}$ for each linking between the $i$-th and $j$-th components. Thus, for a particular choice of sign data, we have the diagram corresponding to Lemma \ref{lem:h} colored by $(-1)^{|\overline{\sigma}^1|^-}\tau_1$ times
	\begin{align*}
		\left(\prod_{i,j} \tau_i^{\|\overline{\sigma}^j\|l_{ij}} \right)
		\left(\prod_{i\neq j} (-1)^{|\overline{\sigma}^i|^-l_{ij}} \right)\nabla(\mathcal{L})((-1)^{|\overline{\sigma}^1|^-}\tau_1 ,\dots, (-1)^{|\overline{\sigma}^n|^-}\tau_n)
		\dfrac{
		(-1)^{|\overline{\sigma}^1|^-}(\tau_1-\tau_1^{-1})}{2(t_{11}-t_{11}^{-1})}.
	\end{align*} By Lemmas \ref{lem:factor} and \ref{lem:h}, this simplifies to
	\begin{align*}
		\left(\prod_{i,j} \tau_i^{\|\overline{\sigma}^j\|l_{ij}} \right)
		\left(\prod_{i\neq j} (-1)^{|\overline{\sigma}^i|^-l_{ij}} \right)
		(-1)^{\sum_i(1+\sum_{j\neq i}l_{ij})|\overline{\sigma}^i|^-}\nabla(\mathcal{L})(\tau_1 ,\dots, \tau_n)\\
		=\left(\prod_{j=1}^n\left(\prod_{i=1}^n \tau_i^{l_{ij}} \right)^{\|\overline{\sigma}^j\|} \right)
		(-1)^{\sum_i|\overline{\sigma}^i|^-}\nabla(\mathcal{L})(\tau_1 ,\dots, \tau_n)
	\end{align*}
	for this choice of sign data. We set $T_j=\prod_{i=1}^n \tau_i^{l_{ij}}$. Thus, 
	\begin{align*}
		\nabla(\mathcal{L}^{\parallel \overline{d}})(t_{1,1},\dots, t_{n,d_n})=\nabla(\mathcal{L})(\tau_1 ,\dots, \tau_n)\sum_{\overline{\sigma}^1, \dots, \overline{\sigma}^n}\left(\prod_j T_j^{^{\|\overline{\sigma}^j\|}}(-1)^{|\overline{\sigma}^j|^-} \right).
	\end{align*}
Note that $\sum_{\overline{\sigma}^1, \dots, \overline{\sigma}^n}\left( \prod_{j=1}^n T_j^{^{\|\overline{\sigma}^j\|}}(-1)^{|\overline{\sigma}^j|^-} \right)$ can be written as $\prod_{j=1}^n\left( \sum_{\overline{\sigma}^j} T_j^{^{\|\overline{\sigma}^j\|}}(-1)^{|\overline{\sigma}^j|^-} \right)$ by expanding and recollecting terms. Indexing the sum over $\overline{\sigma}^j$ by $m_j=|\overline{\sigma}^j|^-$ yields
\begin{align*}
	\prod_{j=1}^n\left( \sum_{\overline{\sigma}^j} T_j^{^{\|\overline{\sigma}^j\|}}(-1)^{|\overline{\sigma}^j|^-} \right)
	=\prod_{j=1}^n\left( \sum_{m_j=0}^{d_j-1}{d_j-1\choose m_j}\left( T_j^{d_j-1-2m_j}(-1)^{m_j} \right) \right)
	=\prod_{j=1}^n (T_j-T_j^{-1})^{d_j-1}.
\end{align*}
After making the necessary substitutions back, the theorem is proven.
\end{proof}

In the case of knots, Theorem \ref{thm:link} can be stated as follows. 

\begin{restatable}{cor}{cable}\label{cor:cable}
	The Conway Potential Function of the $n$-parallel cabling of a knot diagram $\mathcal{K}$ with writhe $w$ is
	\begin{align*}
		\nabla(\mathcal{K}^{\parallel n})(t_1,\dots,t_n)={(t_{1}^w\cdots t_{n}^w-t_{1}^{-w}\cdots t_{n}^{-w})^{n-1}}\nabla(\mathcal{K})(t_{1}\cdots t_{n}).
	\end{align*}
\end{restatable}

\section{Cabled Skein Relation} 
\begin{figure}[h!]
	\includegraphics{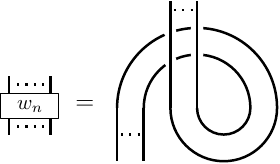}
	\caption{Definition of the full twist on $n$ strands $w_n$.}\label{fig:fulltwist}
\end{figure}

We have a skein relation for cabled strands, which we prove using the diagrammatic calculus. We consider the crossing of cabled strands colored by $t_1,\dots, t_n$ with a writhe adjustment, denoted by $R_{t_1,\dots, t_n}^{\parallel n}$, shown below in Figure \ref{fig:Rpar}. We use a boxed $w_n$, as seen in Figure \ref{fig:fulltwist}, to denote a full twist on $n$ strands. Additional twists are denoted by powers of $w_n$.
\begin{figure}[h!]
	\centering\includegraphics{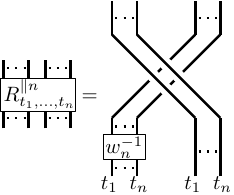}\caption{Definition of the $n$-cabled crossing  $R_{t_1,\dots, t_n}^{\parallel n}$.}\label{fig:Rpar}
\end{figure}
\begin{lem}\label{lem:twist}
	Fix sign data $\overline{\sigma}\in\{\pm1\}^{n-1}$. The endomorphism $\bracks{w_n}$ acts on the image of $\v_{\overline{\sigma}}^{t_1,\dots, t_n}$ in $\bigotimes_{i=1}^n V(t_{i})$ by the scalar $\tau^{\|\overline{\sigma}\|}$.
\end{lem}
\begin{proof}
	The proof follows immediately from Lemma \ref{lem:slide} and the relations in Figure \ref{fig:maxmin}. 
\end{proof}
\cableskein
\begin{proof}
	Our approach will once again make use of the diagrammatic calculus developed here. As in the proof of Theorem \ref{thm}, we will verify the relation for any choice of sign data on our fork diagram.\newline
	
	Set $\tau=t_1\cdots t_n$. Select sign data $\overline{\sigma}^1,\overline{\sigma}^2\in\{\pm 1\}^{n-1}$, and $\sigma^3\in\{\pm 1\}$, and consider the fork \[\left(\v_{\overline{\sigma}^1}^{t_1,\dots, t_n}\otimes \v_{\overline{\sigma}^2}^{t_1,\dots, t_n}\right)\circ \v_{{\sigma^3}}^{(-1)^{|\overline{\sigma}^1|^-}\tau,(-1)^{|\overline{\sigma}^2|^-}\tau}:V\left((-1)^{|\overline{\sigma}^1|^-+|\overline{\sigma}^2|^-+|{\sigma}^3|^-}\tau^2\right)\to \bigotimes_{i=1}^n V(t_i)\otimes \bigotimes_{i=1}^n V(t_i).\]
	Making use of Lemma \ref{lem:twist} and Remark \ref{rem:slidelink}, the composition of a fork with $\left(R_{t_1,\dots, t_n}^{\parallel n}\right)^2$ yields $\left(t_1^2\cdots t_n^2\right)^{2\sigma^3}$, as shown in Figure \ref{fig:Rfork} below. Since we obtain the inverse factor on $\left(R_{t_1,\dots, t_n}^{\parallel n}\right)^{-2}$, the skein relation is valid when composed with a fork for any sign data. Thus, the equality is true in general.
	\begin{figure}[h!]
		\centering\includegraphics{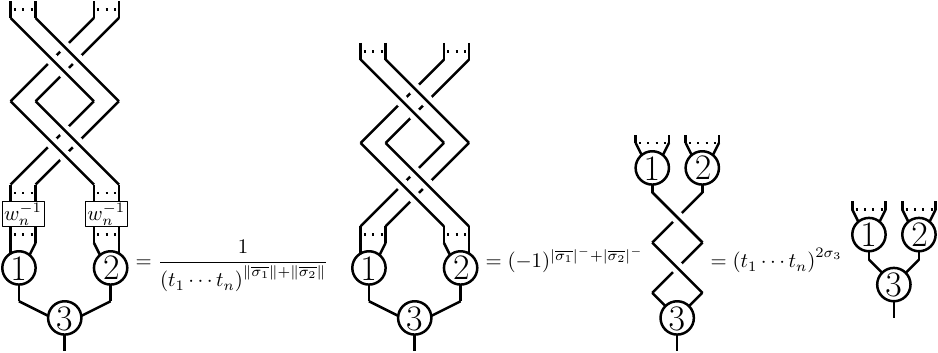}\caption{Action of $\left(R_{t_1,\dots, t_n}^{\parallel n}\right)^2$ on a fork.}\label{fig:Rfork}
	\end{figure}
	
\end{proof}
\begin{ex}
	Using this relation, we compute the CPF of an $n$-cabled trefoil with writhe $w$. Note that the trace of ${\left(R_{t_1,\dots, t_n}^{\parallel n}\right)^\pm (\bracks{w_n}^w\otimes id)}$ corresponds to a cabled unknot twisted $w$ times. The CPF of these unknots equals $\dfrac{\left(\tau^{w}-\tau^{-w}\right)^{n-1}}{\tau-\tau^{-1}}$, having set $\tau=t_1\cdots t_n.$ Then
\begin{align}
	\left(R_{t_1,\dots, t_n}^{\parallel n}\right)^3(\bracks{w_n}^w\otimes id)&=
\left(
\left(\tau^2+\tau^{-2}\right)\left(R_{t_1,\dots, t_n}^{\parallel n}\right)-\left(R_{t_1,\dots, t_n}^{\parallel n}\right)^{-1}\right)(\bracks{w_n}^w\otimes id),\\
\nabla(\mathsf{3_1}^{\parallel n})(t_1,\dots, t_n)&=\left(\tau^w-\tau^{-w}\right)^{n-1}
\dfrac{\left(\tau^2-1+\tau^{-2}\right)}{\tau-\tau^{-1}}=\left(\tau^w-\tau^{-w}\right)^{n-1}\nabla(\mathsf{3_1})(\tau),
\end{align}
and this result is consistent with Corollary \ref{cor:cable}.
\end{ex}

\bibliographystyle{alpha}

\end{document}